\documentclass[oneside,11pt,reqno]{amsart}

\usepackage{verbatim,upref,amsxtra,amssymb,amscd}

\usepackage{varioref}

\DeclareMathAlphabet\mathscr{U}{eus}{m}{n}
\SetMathAlphabet\mathscr{bold}{U}{eus}{b}{n}
\DeclareMathAlphabet\matheur{U}{eur}{m}{n}
\SetMathAlphabet\matheur{bold}{U}{eur}{b}{n}
\newcommand{\IGNORE}[1]{} 
\newcommand{\Z}{\mathbb{Z}} \newcommand{\R}{\mathbb{R}}

  \newcommand{\N}{\mathbb{N}}
  \newcommand {\absolute}[1] {\left| {#1} \right|}
\newcommand {\norm}[1] {\left\| {#1} \right\|}

\numberwithin{equation}{section}

\newtheorem{theo}{Theorem}[section]
\newtheorem{prop}[theo]{Proposition} \newtheorem{lemm}[theo]{Lemma}

  \newtheorem{coro}[theo]{Corollary}

\theoremstyle{definition}

\newtheorem{defi}[theo]{Definition} 
\newtheorem{exas}[theo]{Examples}

\theoremstyle{remark}

\newtheorem{rema}[theo]{Remark}

\begin{document}
\allowdisplaybreaks\frenchspacing \setlength{\baselineskip}{1.1\baselineskip}


\title[Invariant sets and measures]{Invariant sets and measures of
nonexpansive group automorphisms}

\author{Elon Lindenstrauss}

\address{Elon Lindenstrauss: Department of Mathematics, Stanford
University,\linebreak Stanford, CA 94305, USA}
\email{elonl@math.stanford.edu}

\author{Klaus Schmidt}

\address{Klaus Schmidt: Mathematics Institute, University of Vienna,
Strudlhofgasse~4, A-1090 Vienna, Austria \newline\indent \textup{and}
\newline\indent Erwin Schr{\"o}dinger Institute for Mathematical
Physics, Boltzmanngasse~9, A-1090 Vienna, Austria}
\email{klaus.schmidt@univie.ac.at}


    \begin{abstract}
We prove that the restriction of a probability measure invariant
under a nonhyperbolic, ergodic and  totally irreducible automorphism
of a compact connected abelian group to the leaves of the central
foliation is severely restricted. We also prove a topological
analogue of this result: the intersection of every proper closed
invariant subset with each central leaf is compact.
    \end{abstract}

\maketitle

\section{Introduction}
    \label{s:intro}

A continuous automorphism $\alpha $ of an additive compact abelian
group $X$ is \emph{expansive} if there exists a neighbourhood $N(0)$
of the identity $0\in X$ with $\bigcap_{n\in\mathbb{Z}}\alpha
^n(N(0))=\{0\}$, \emph{irreducible} if every closed $\alpha
$-invariant subgroup $Y\subsetneq X$ is finite, \emph{totally
irreducible} if every nonzero power of $\alpha $ is irreducible, and
\emph{ergodic} if it is topologically transitive (and hence ergodic
with respect to the normalized Haar measure $\lambda _X$ of $X$).

In this paper we study the collection of invariant measures of a
nonhyperbolic, ergodic and totally irreducible automorphism of the
$n$-torus $\mathbb{T}^n$ or, more generally, of a nonexpansive,
ergodic and totally irreducible automorphism $\alpha $ of a compact
connected abelian group $X$. Every nonhyperbolic, ergodic and
irreducible automorphism $\alpha $ of $\mathbb{T}^n$ is partially
hyperbolic in the usual sense\footnote{A $C^2$ diffeomorphism $f$ of
a Riemannian manifold $M$ is \emph{partially hyperbolic} if there is
a $Df$-invariant splitting $TM=E_s \oplus  E_c \oplus E_u$ of the
tangent manifold $TM$ of $M$ in which at least two of the sub-bundles
are nontrivial, so that $Df$ uniformly expands all vectors in $E_u$,
uniformly contracts all vectors in $E_s$, and the vectors in $E_c$
are neither expanded as strongly as any vector in $E_u$ nor
contracted as strongly as any vector in $E_s$.} with the additional
property that its derivative $D\alpha$ preserves the length of
vectors in $E_c$. For an arbitrary nonexpansive, ergodic and totally
irreducible continuous automorphisms $\alpha $ of a compact connected
abelian group $X$ these `sub-bundles' can be more complicated objects
(due to the fact that the group need not be locally connected), but
an analogue of this strong form of partial hyperbolicity also holds
in this more general situation.

\IGNORE{Surprisingly, and in contrast to the hyperbolic case, very
little is known about nonhyperbolic toral automorphisms. For example,
despite the recent progress in understanding ergodicity of partially
hyperbolic maps (see \cite {PS} and \cite {BPSW}) virtually nothing
is known about the dynamics of a perturbation of an irreducible
nonhyperbolic toral automorphism.

We investigate a different aspect of the dynamics of this
prototypical nonhyperbolic system, namely trying to understand the
invariant measures on this system.}

Let $\alpha $ be a nonexpansive, ergodic and totally irreducible
automorphism of a compact connected abelian group $X$. The normalized
Haar measure $\lambda _X$ of $X$ is obviously invariant under $\alpha
$, and Y. Katznelson \cite{Katznel} proved that the
measure-preserving system $(X,\alpha,\mathscr{B}_X,\lambda_X)$ (where
$\mathscr{B}_X$ denotes the Borel sigma-algebra of $X$) is
measure-theoretically isomorphic to a Bernoulli shift.

There is another family of --- admittedly not very interesting ---
$\alpha $-invari\-ant ergodic probability measures on $X$: let
$X^{(0)}\subset X$ be the dense central subgroup of $\alpha $ defined
in \eqref{eq:Sc0}, on which $\alpha $ acts isometrically. Then the
closure of the $\alpha $-orbit of any element $x\in X^{(0)}$ is a
compact $\alpha $-invariant subset of $X^{(0)}$ (and hence of $X$) on
which $\alpha $ acts with a unique $\alpha $-invariant measure
denoted by $\tilde{\lambda }_x$.

It is not immediate how to construct other invariant measures; in
fact, the main result in this paper shows that \emph{all} $\alpha
$-invariant probability measures $\mu \ne \lambda _X$ on $X$ satisfy
a somewhat surprising rigidity phenomenon related to the scarcity of
invariant measures under a multidimensional abelian semigroup of
toral endomorphisms. This scarcity of invariant measures was
conjectured by H. Furstenberg and is still open, though there are
important partial results by several authors including D. Rudolph
\cite{Rudolph} for the one-dimensional case and A. Katok and R.
Spatzier \cite{KS} in the higher-dimensional case.

In order to describe this rigidity property we use a construction
from \cite{KS} to define a system of `conditional' measures on the
leaves of the central foliation induced by an $\alpha$-invariant
measure $\mu $ on $X$. In general, if we start with an $\alpha
$-invariant probability measure $\mu $ on $X$, these leaf measures
will only be sigma-finite. Indeed, for $\mu =\lambda _X$, the induced
measure on each central leaf is the (infinite) Haar measure on the
leaf. Our main result is that the leaf measures are finite for any
$\alpha $-invariant probability measure $\mu $ on $X$ which does not
contain a copy of $\lambda _X$ in its ergodic decomposition.

\medskip\noindent \textbf{Theorem} (Theorem \ref{t:finite}).\,
\emph{Let $\alpha$ be a nonexpansive, ergodic and totally irreducible
automorphism of a compact connected abelian group $X$ with normalized
Haar measure $\lambda _X$, and let $\mu $ be an $\alpha $-invariant
probability measure on $X$ which is singular with respect to $\lambda
_X$. Then the conditional measure $\rho _x$ on the central leaf
through $x$ \textup{(}defined in \eqref{eq:rhox}\textup{)} is finite
for almost every $x\in X$.}

Both the statement and the proof of Theorem \ref{t:finite} are
modelled on Host's proof of Rudolph's Theorem in \cite{Host1} and its
generalization in \cite{Host2}.

\medskip The following two definitions can easily be adapted to the
general setting of partially hyperbolic maps.

    \begin{defi}
    \label{d:equivalent meas}
Two $\alpha$-invariant probability measures $\mu_1,\mu_2$ on $X$ are
\emph{centrally equivalent} if they have an invariant joining $\nu $
(i.e. an $(\alpha\times\alpha)$-invariant measure $\nu $ on $X\times
X$ which projects to $\mu_1$ and $\mu_2$, respectively) so that, for
$\nu\textsl{-a.e.}\;(x,y)\in X\times X$, $x$ and $y$ lie on the same
central leaf; in other words,
    $$
x-y\in X^{(0)}\enspace \textup{for}\enspace \nu
\textsl{-a.e.}\;(x,y)\in X\times X,
    $$
where $X^{(0)}\subset X$ is the central subgroup of $\alpha $ defined
in \eqref{eq:Sc0}.
    \end{defi}

    \begin{defi}
An $\alpha$ invariant probability measure $\mu $ on $X$ is
\emph{virtually hyperbolic} if there exists an $\alpha $-invariant
Borel set $Z\subset X$ with $\mu (Z)=1$ which intersects every
central leaf in at most one point, i.e. with $Z\cap(x+Z)=\varnothing
$ for every $x\in X^{(0)}$.
    \end{defi}

In Section \ref{s:virtually hyperbolic} we prove that Theorem
\ref{t:finite} implies the following result.

    \begin{theo}
    \label{t:measures}
Let $\alpha $ be a nonexpansive, ergodic and totally irreducible
automorphism of a compact connected abelian group $X$ with normalized
Haar measure $\lambda _X$, and let $\mu $ be an $\alpha $-invariant
probability measure on $X$ which is singular with respect to $\lambda
_X$. Then the following conditions are satisfied.
    \begin{enumerate}
    \item
There is a virtually hyperbolic $\alpha $-invariant probability
measure $\mu '$ on $X$ which is centrally equivalent to $\mu $;
    \item
If $\mu $ is weakly mixing \textup{(}or, more generally, if the point
spectrum of the action of $\alpha $ on $L^2 (X,\mathscr{S}, \mu )$
contains no eigenvalue of $\alpha $ of absolute value $1$\textup{)},
then $\mu$ is virtually hyperbolic;
    \item
If $\mu$ is ergodic, but not necessarily weakly mixing, we write, for
every $x\in X^{(0)}$, $\tilde{\lambda }_x$ for the unique $\alpha
$-invariant probability measure on $X^{(0)}$ --- and hence on $X$ ---
concentrated on the compact orbit closure $\overline{\{\alpha
^nx:n\in \mathbb{Z}\}}$ of $x$ under $\alpha $. Then $\mu $ is an
ergodic component of $\mu ' * \tilde{\lambda }_{x_0}$ for some
$x_0\in X^{(0)}$.
    \end{enumerate}
    \end{theo}

Finally, in Section \ref{s:topological} we prove the following
topological analogue of the Theorems \ref{t:measures} and
\ref{t:finite}.

\medskip\noindent \textbf{Theorem} (Theorem \ref{theorem:
topological}).\, \emph{Let $\alpha$ be a nonexpansive, ergodic and
totally irreducible automorphism of a compact connected abelian group
$X$. Then any closed $\alpha$-invariant subset $Y \subsetneq X$
intersects every central leaf in a compact subset of the leaf.}

\subsection*{Acknowledgements}
This research has been supported in part by NSF grant DMS 0140497
(E.L.) and FWF Project P14379-MAT (K.S.). During part of this work,
both authors received support from the American Institute of
Mathematics and NSF grant DMS 0222452.  We would furthermore like to
express our gratitude to the Mathematics Department, University of
Washington, Seattle, the Newton Institute, Cambridge and ETH Z\"urich
for hospitality during parts of this work, and Doug Lind for some
interesting and useful discussions.

\section{Irreducible group automorphisms}
    \label{s:irreducible}

Let $\alpha $ and $\beta $ be continuous automorphisms of compact
abelian groups $X$ and $Y$, respectively. Then $\alpha $ and $\beta $
are \emph{conjugate} if there exists a continuous group isomorphism
$\phi \colon X\longrightarrow Y$ with
    \begin{equation}
    \label{eq:conjugate}
\beta \circ \phi =\phi \circ \alpha ,
    \end{equation}
and $\beta $ is a \emph{factor} of $\alpha $ if there exists a
continuous surjective group homomorphism $\phi \colon
X\longrightarrow Y$ satisfying \eqref{eq:conjugate}. The map $\phi $
in \eqref{eq:conjugate} is called an (\emph{algebraic})
\emph{conjugacy} or an (\emph{algebraic}) \emph{factor map}. The
automorphisms $\alpha $ and $\beta $ are \emph{weakly conjugate} if
each of them is a factor of the other, and \emph{finitely equivalent}
if each of them is a factor of the other with a finite-to-one factor
map.

We recall a few basic facts about irreducible ergodic automorphisms
of compact abelian groups. Let $R_1=\mathbb{Z}[u^{\pm1}]$ be the ring
of Laurent polynomials with integral coefficients. We write every
$h\in R_1$ as
    \begin{equation}
    \label{eq:h}
\smash[t]{h=\sum_{m\in\mathbb{Z}}h_mu^m}
    \end{equation}
with $h_m\in\mathbb{Z}$ for every $m\in\mathbb{Z}$ and $h_m=0$ for
all but finitely many $m$.

Let $\alpha $ be an automorphism (always assumed to be continuous) of
a compact abelian group $X$ with (additive) dual group $\hat{X}$, and
let $\hat{\alpha }$ be the dual automorphism of $\hat{X}$ defined by
    $$
\langle \hat{\alpha }a,x\rangle =\langle a,\alpha x\rangle
    $$
for every $x\in X$ and $a\in \hat{X}$, where $\langle a,x\rangle $
denotes the value of $a\in\hat{X}$ at $x\in X$. For every
$h=\sum_{n\in\mathbb{Z}} h_nu^n\in R_1$, $x\in X$ and $a\in\hat{X}$
we set
    \begin{equation}
    \label{eq:h(alpha)}
\smash[b]{h(\alpha )(x)=\sum_{n\in\mathbb{Z}} h_n\alpha ^nx,\qquad
h(\hat{\alpha })(a)=\sum_{n\in\mathbb{Z}} h_n\hat{\alpha }^na,}
    \end{equation}
and note that
    \begin{equation}
    \label{eq:m0}
\langle h(\hat{\alpha })(a),x\rangle =\langle \widehat{h(\alpha
)}(a),x\rangle =\langle a,h(\alpha )(x)\rangle.
    \end{equation}
The dual group $\hat{X}$ is a module over the ring $R_1$ with operation
    \begin{equation}
    \label{eq:m1}
h\cdot a=h(\hat{\alpha })(a)
    \end{equation}
for $h\in R_1$ and $a\in\hat{X}$. In particular,
    \begin{equation}
    \label{eq:m2}
u^m\cdot a = \hat{\alpha }^ma
    \end{equation}
for $m\in\mathbb{Z}$ and $a\in \hat{X}$. This module is called the
\emph{dual module} $M=\hat{X}$ of $\alpha $. Conversely, if $M$ is an
$R_1$-module, we obtain an automorphism $\alpha _M$ on the compact
abelian group
    \begin{equation}
    \label{eq:m3}
X_M=\widehat{M}
    \end{equation}
whose dual automorphism is defined by
    \begin{equation}
    \label{eq:m4}
\hat{\alpha }_Ma=u\cdot a
    \end{equation}
for every $a\in M$.

    \begin{exas}[\cite{DSAO}]
    \label{e:autos}
(1) Let $M=R_1$. Since $R_1$ is isomorphic to the direct sum
$\sum_{\mathbb{Z}}\mathbb{Z}$ of copies of $\mathbb{Z}$, indexed by
$\mathbb{Z}$, the dual group $X=\widehat{R_1}$ is isomorphic to the
cartesian product $\mathbb{T}^{\mathbb{Z}}$ of copies of
$\mathbb{T}=\mathbb{R}/\mathbb{Z}$. We write a typical element
$x\in\mathbb{T}^{\mathbb{Z}}$ as $x=(x_{n})$ with
$x_{n}\in\mathbb{T}$ for every $n\in\mathbb{Z}$ and choose the
following identification of $X_{R_1}=\widehat{R_1}$ and
$\mathbb{T}^{\mathbb{Z}}$: for every $x=(x_{n})$ in
$\mathbb{T}^{\mathbb{Z}}$ and $h=\sum_{n\in\mathbb{Z}}h_nu^n\in R_1$,
    \begin{equation}
    \label{eq:hatR}
\langle x,h\rangle=e^{2\pi i \sum_{n\in\mathbb{Z}} h_n x_{n}}.
    \end{equation}
Under this identification the automorphism $\alpha _{R_1}$ on
$X_{R_1}=\mathbb{T}^{\mathbb{Z}}$ becomes the shift
    \begin{equation}
    \label{eq:tau}
(\tau x)_m=x_{m+1}
    \end{equation}
with $m\in\mathbb{Z}$ and $x=(x_m)\in X_{R_1}= \mathbb{T}^{\mathbb{Z}}$.

\smallskip (2) Let $I\subset R_1$ be an ideal, and let $M=R_1/I$.
Since $M$ is a quotient of the additive group $R_1$ by an $\hat\alpha
_{R_1}$-invariant subgroup, the dual group $X_M$ is the $\alpha
_{R_1}$-invariant subgroup
    \begin{equation}
    \label{eq:ideal}
    \begin{aligned}
X_{R_1/I}&=I^\perp=\{x\in X_{R_1}=\mathbb{T}^{\mathbb{Z}}: \langle
x,h\rangle=1\enspace\text{for every}\enspace h\in I\}
    \\
&=\smash[b]{\biggl\{x\in\mathbb{T}^{\mathbb{Z}}:
\sum_{n\in\mathbb{Z}} h_nx_{m+n}=0\pmod1}
    \\
&\qquad\qquad\qquad\qquad\qquad\qquad\smash{\text{for every}\enspace
h\in I\enspace\text{and}\enspace m\in\mathbb{Z}\biggr\}}
    \\
&=\{x\in\mathbb{T}^{\mathbb{Z}}: h(\tau )(x)=0\enspace \text{for
every}\enspace h\in I\},
    \end{aligned}
    \end{equation}
and $\alpha _{R_1/I}$ is the restriction of $\tau =\alpha _{R_1}$ to
$X_{R_1/I}\subset\mathbb{T}^{\mathbb{Z}}=X_{R_1}$.

We can express \eqref{eq:ideal} as
    \begin{equation}
    \label{eq:kernels}
X_{R_1/I}=\widehat{X/I}=I^\perp=\bigcap_{h\in I}\ker(h(\tau )).
    \end{equation}

If $I=(f)=fR_1$ is the principal ideal generated by some $f\in R_1$,
then \eqref{eq:kernels} becomes
    \begin{equation}
    \label{eq:kernel}
X_{R_1/(f)}=\widehat{X/(f)}=(f)^\perp=\ker(f(\tau )).
    \end{equation}

\smallskip (3) Let $\alpha $ be the automorphism of the $m$-torus
$X=\mathbb{T}^m=\mathbb{R}^m/\mathbb{Z}^m$ defined by a matrix
$A\in\textup{GL}(m,\mathbb{Z})$. Then the dual module $M=\hat{X}$ is
equal to $\mathbb{Z}^m$ with operation $f\cdot
\mathbf{m}=f(A^\top)(\mathbf{m})$ for every $f\in R_1$ and
$\mathbf{m}\in\mathbb{Z}^m$ (cf. \eqref{eq:m1}), where
$A^\top\in\textup{GL}(m,\mathbb{Z})$ is the transpose matrix of $A$.

The automorphism $\alpha $ is irreducible if and only if the
characteristic polynomial $f=f_0+\dots +f_{m-1}u^{m-1}+u^m$ of $A$ is
irreducible, and $\alpha $ is conjugate to $\alpha _{R_1/(f)}$ if and
only if $A$ is conjugate in $\textup{GL}(m,\mathbb{Z})$ to the
companion matrix
    \begin{equation}
    \label{eq:companion}
C_f=\left(
    \begin{smallmatrix}
0&1&\cdots&0&0
    \\
0&0&\cdots&0&0
    \\
\vdots&\vdots&\ddots&\vdots&\vdots
    \\
0&0&\cdots&0&1
    \\
-f_0&-f_1&\hdots&-f_{m-2}&-f_{m-1}
    \end{smallmatrix}
\right)\in\textup{GL}(m,\mathbb{Z}).
    \end{equation}
    \end{exas}

    \begin{theo}
    \label{t:irreducible}
Let $\alpha $ be an irreducible automorphism of an infinite compact
connected abelian group $X$. Then there exists a unique irreducible
polynomial $f=f_0+\dots +f_nu^n\in R_1$ with the following properties.
    \begin{enumerate}
    \item
$n\ge1$, $f_n>0$ and $f_0\ne0$;
    \item
$\alpha $ is finitely equivalent to $\alpha _{R_1/(f)}$, where
$(f)=fR_1\subset R_1$ is the ideal generated by $f$ \textup{(}cf.
Example \ref{e:autos} \textup{(2))};
    \item
$\alpha $ is ergodic if and only if $f$ is not cyclotomic
\textup{(}i.e. $f$ does not divide $u^m-1$ for any $m\ge1$\textup{)};
    \item
$\alpha $ is expansive if and only if $f$ has no roots of absolute value $1$.
    \item
$\alpha$ is totally irreducible if and only if $f$ has no two
distinct roots whose ratio is a root of unity.
    \end{enumerate}
Conversely, if $f=f_0+\dots +f_nu^n\in R_1$ is an irreducible
polynomial satisfying condition \textup{(1)} above, then the group
$X_{R_1/(f)}$ in \eqref{eq:ideal} is connected and the automorphism
$\alpha _{R_1/(f)}$ of $X_{R_1/(f)}$ is irreducible.
    \end{theo}

    \begin{proof}
The statements (1)--(4) and the converse follow from
\cite[Proposition 2.7 and Theorem 29.2]{DSAO}.

For the proof of (5) we note that $\alpha $ is totally irreducible if
and only if $\alpha _{R_1/(f)}$ is totally irreducible. Since $\alpha
_{R_1/(f)}^m$ is dual to multiplication by $u^m$ on
$\hat{X}=R_1/(f)=M$, $\alpha _{R_1/(f)}^m$ is irreducible if and only
if the subgroup
    $$
N=\{h(u^m):h\in R_1\}/(f)\subset R_1/(f)=M
    $$
has finite index in $R_1/(f)$. As the group $M$ is torsion-free, the
latter condition is equivalent to the statement that
$N\otimes_\mathbb{Z}\mathbb{Q}=
M\otimes_\mathbb{Z}\mathbb{Q}\cong\mathbb{Q}^n$, and hence to the
condition that the elements $\{(1+(f)),(u^m+(f)),\dots
,(u^{m(n-1)}+(f))\}$ in $M$ are rationally independent. In other
words, $\alpha _{R_1/(f)}^m$ is reducible if and only if one can find
a nonzero element $(k_0,\dots ,k_{n-1})\in\mathbb{Z}^n$ with
    \begin{equation}
    \label{eq:dependence}
g(u^m)=k_0+k_1u^m+\dots +k_{n-1}u^{m(n-1)}\in(f),
    \end{equation}
where we may assume without loss of generality that the resulting
polynomial $g\in R_1$ is irreducible. By evaluating
\eqref{eq:dependence} on any root $\theta $ of $f$ we obtain that
$g(\theta ^m)=0$ for every root $\theta $ of $f$, and Galois theory
shows that the degree of $g$ is equal to the number of distinct
elements in the set $\Omega _f^{(m)}=\{\theta ^m:\theta \enspace
\textup{is a root of}\enspace f\}$. This proves that $\alpha
_{R_1/(f)}^m$ is irreducible if and only if the cardinality of
$\Omega _f^{(m)}$ is equal to $n$, which implies (5).
    \end{proof}

Example \ref{e:autos} (2) gives an explicit representation --- up to
finite equivalence --- of every irreducible automorphism of a compact
connected abelian group $X$. For an alternative description we follow
\cite{ES} (for background see \cite{S2}, \cite[Section 7]{DSAO} and
\cite{Weil}).

Let $\alpha $ be an irreducible automorphism of an infinite compact
connected abelian group $X$, and let $f\in R_1$ be the irreducible
polynomial appearing in Theorem \ref{t:irreducible}. We fix a root
$\theta \in\bar{\mathbb{Q}}$ of $f$, denote by $K=\mathbb{Q}(\theta
)$ the algebraic number field generated by $\theta $, and write
$P^{(K)}$, $P^{(K)}_f$, and $P^{(K)}_\infty$, for the sets of places
(= equivalence classes of valuations), finite places and infinite
places of $K$. For every place $v$ of $K$ and every valuation $\phi
\in v$, the \emph{$v$-adic completion} $K_v$ of $K$ (i.e. the
completion of $K$ with respect to metric $\delta (a,b)=\phi (a-b)^\gamma
$ for some suitable $\gamma >0$) is a locally compact, metrizable
field and hence a locally compact additive group. We fix a Haar
measure $\lambda_v$ on the additive group $K_v$ and denote by
$\textup{mod}_{K_v}\colon K_v\longrightarrow \mathbb{R}$ the map
satisfying
    \begin{equation}
    \label{eq:valuation}
\lambda _v(aB)=\textup{mod}_{K_v}(a)\lambda _v(B)
    \end{equation}
for every $a\in K_v$ and every Borel set $B\subset K_v$. The
restriction of $\textup{mod}_{K_v}$ to $K$ is a valuation in $v$,
denoted by $|\cdot |_v$.

Let
    \begin{equation}
    \label{eq:Pc}
P=\{v\in P^{(K)}_{f}:|\theta |_v\ne1\},\enspace S=P^{(K)}_\infty \cup P.
    \end{equation}
For every infinite place $v\in P_\infty^{(K)}$, the $v$-adic
completion $K_v$ is either equal to $\mathbb{R}$ or to $\mathbb{C}$
(in particular, $K_v=\mathbb{C}$ for any $v\in S^{(0)}$). We write
    \begin{equation}
    \label{eq:iotav}
\iota _v\colon K\longrightarrow K_v (= \mathbb{R}\enspace
\textup{or}\enspace \mathbb{C})
    \end{equation}
for the embedding of $K$ in its completion $K_v$ and use the same
symbol $\iota _v$ to denote the corresponding identification of $K_v$
with $\mathbb{R}$ or $\mathbb{C}$.

The set
    \begin{equation}
    \label{eq:W}
\smash[t]{W=\prod_{v\in S}K_v}
    \end{equation}
is a locally compact algebra over $K$ with respect to coordinate-wise
addition, multiplication and scalar multiplication (with scalars in
$K$). We write every $w\in W$ as $w=(w_v)=(w_v,\,v\in S)$ with
$w_v\in K_v$ for every $v\in S$ and define
    \begin{equation}
\label {eq:norm} \left\| w \right\| =\max_{v\in S}| w_v |_v.
    \end{equation}
Let $\bar{\beta }$ be the automorphism of $W$ given by
    \begin{equation}
    \label{eq:barbetac}
\bar{\beta }w=(\theta w_v)
    \end{equation}
for every $w=(w_v)\in W$.

We put
    \begin{equation}
    \label{eq:Rc}
\mathscr{R}=\{a\in K:|a|_v\le1\enspace\text{for every}\enspace v\in
P^{(K)}\smallsetminus S\}\supset\mathfrak{o}_{K},
    \end{equation}
where $\mathfrak{o}_K$ is the ring of integers in $K$, and denote by
    \begin{equation}
    \label{eq:ic}
\iota \colon K\longrightarrow W
    \end{equation}
the diagonal embedding $a\mapsto \iota (a)=(a,\dots,a)$, $a\in K$. By
abuse of notation we identify each $K_v,\,v\in S$, with the subgroup
    $$
\{w\in W: w_{v'} =0 \text{ for every }v'\neq v\}\subset W.
    $$

    \begin{theo}
    \label{t:irreducible2}
Suppose that $\alpha $ is an automorphism of an infinite compact
connected abelian group $X$. Then $\alpha $ is irreducible if and
only if there exist an element $\theta \in\bar{\mathbb{Q}}^\times
=\bar{\mathbb{Q}}\smallsetminus \{0\}$ and a finitely generated
$\mathbb{Z}[\theta ^{\pm1}]$-submodule $L\subset K=\mathbb{Q}(\theta
)$ such that $\alpha $ is algebraically conjugate to the automorphism
$\beta _{(\theta ,L)}$ on the quotient group
    \begin{equation}
    \label{eq:YL}
Y_L=W/\iota (L)
    \end{equation}
induced by $\bar{\beta }$ \textup{(}cf.
\eqref{eq:Pc}--\eqref{eq:barbetac}\textup{)}.

    \begin{enumerate}
    \item
The following conditions are equivalent.
    \begin{enumerate}
    \item
$\alpha $ is ergodic,
    \item
$\theta $ is not a root of unity.
    \end{enumerate}

\smallskip
    \item
The following conditions are equivalent.
    \begin{enumerate}
    \item
$\alpha $ is expansive,
    \item
The orbit of $\theta $ under the action of the Galois group
$\textup{Gal}[\bar{\mathbb{Q}}:\mathbb{Q}]$ does not intersect
$\mathbb{S}=\{z\in\mathbb{C}:|z|=1\}$.
    \end{enumerate}

\smallskip
    \item
The following conditions are equivalent.
    \begin{enumerate}
    \item
$X\cong \mathbb{T}^n$ for some $n\ge1$,
    \item
$S=P_\infty ^{(K)}$,
    \item
$\theta $ is an algebraic unit.
    \end{enumerate}
    \item
The following conditions are equivalent.
    \begin{enumerate}
    \item
$\alpha$ is totally irreducible,
    \item
The orbit of $\theta $ under the action of the Galois group
$\textup{Gal}[\bar{\mathbb{Q}}:\mathbb{Q}]$ does not contain two
distinct elements whose ratio is a root of unity.
    \end{enumerate}
    \end{enumerate}
    \end{theo}

    \begin{proof}
\cite[Corollary 3.5]{ES}, \cite[Theorem 7.1 and Propositions
7.2--7.3]{DSAO} and Theorem \ref{t:irreducible} in this paper.
    \end{proof}

    \begin{rema}
    \label{r:unit}
If $\theta $ is an algebraic unit, then $S=P_\infty ^{(K)}$ and
    \begin{equation}
    \label{eq:infty}
W\cong\mathbb{R}^{r(K)},\qquad Y\cong\mathbb{T}^{r(K)},
    \end{equation}
where
    \begin{equation}
    \label{eq:r(K)}
r(K)=|\{v\in P_\infty ^{(K)}:K_v=\mathbb{R}\}|+2|\{v\in P_\infty
^{(K)}:K_v=\mathbb{C}\}|.
    \end{equation}

Conversely, if $Y\cong\mathbb{T}^m$ for some $m\ge1$, then $\theta $
is an algebraic unit.
    \end{rema}

\section{Structure and examples of nonexpansive automorphisms}
    \label{s:measures}

Let $\alpha $ be a nonexpansive irreducible ergodic automorphism of a
compact connected abelian group $X$. We apply the Theorem
\ref{t:irreducible2} and assume that
    \begin{equation}
    \label{eq:notation}
\alpha =\beta _{(\theta ,L)},\qquad X=W/\iota (L),
    \end{equation}
for some $\theta \in\bar{\mathbb{Q}^\times }$ and some finitely
generated $\mathbb{Z}[\theta ^{\pm1}]$-submodule $L\subset
K=\mathbb{Q}(\theta )$. Denote by $\lambda _X$ the normalized Haar
measure of $X$ and write
    \begin{equation}
    \label{eq:pi}
\pi \colon W\longrightarrow X=W/\iota (L)
    \end{equation}
for the quotient map (cf. \eqref{eq:Pc}--\eqref{eq:YL}). In the
notation of \eqref{eq:Pc} and \eqref{eq:W} we set
    \begin{equation}
    \label{eq:Sc0}
    \begin{aligned}
S^{(0)}&=\{v\in S:|\theta |_v=1\}\subset P_\infty ^{(K)},
    \\
W^{(0)}&=\{w=(w_v)\in W:w_v=0\enspace \textup{for every}\enspace v\in
S\smallsetminus S^{(0)}\}
    \\
&\cong \prod_{v\in S^{(0)}}K_v\cong \mathbb{C}^{|S^{(0)}|},
    \\
\IGNORE {S^{(u)} &= \{v\in S:|\theta |_v>1\}
    \\
W^{(u)}&=\{w=(w_v)\in W:w_v=0\enspace \textup{for every}\enspace v\in
S\smallsetminus S^{(u)}\}
    \\
S^{(s)} &= \{v\in S:|\theta |_v<1\}
    \\
W^{(s)}&=\{w=(w_v)\in W:w_v=0\enspace \textup{for every}\enspace v\in
S\smallsetminus S^{(s)}\}
    \\
W^*&=\{w=(w_v)\in W:w_v=0\enspace \textup{for every}\enspace v\in S^{(0)}\}
    \\
&\cong \prod_{v\in S\smallsetminus S^{(0)}}K_v,
    \\
} X^{(0)}&=\pi (W^{(0)}).
    \end{aligned}
    \end{equation}

The \emph{central subgroup} group $X^{(0)}\subset X$ is $\alpha
$-invariant and dense by irreducibility. Furthermore, since $|L/\iota
(\mathscr{R})|<\infty $ (cf. \eqref{eq:Rc} and \cite{ES}) and $\iota
(\mathscr{R})\cap W^{(0)}=\{0\}$ by the product formula
(\cite[Theorem 10.2.1]{Cassels}), $L\cap W^{(0)}=\{0\}$.

    \begin{exas}
    \label{e:central}
(1) Let $\alpha $ be a nonexpansive irreducible ergodic automorphism
of $X=\mathbb{T}^m$ defined by a matrix
$A\in\textup{GL}(m,\mathbb{Z})$ with real eigenvalues $\theta
_1,\dots ,\theta _{m_1}$ and complex eigenvalues $\theta
_{m_1+1},\bar{\theta }_{m_1+1},\dots ,\theta _{m_1+m_2},\bar{\theta
}_{m_1+m_2}$, where $m=m_1+2m_2$, and where $\bar{\theta }_i$ is the
complex conjugate of $\theta _i$ for $i=m_1+1,\dots ,m_1+m_2$. We fix
an eigenvalue $\theta $ of $A$, set $K=\mathbb{Q}(\theta )$, and
obtain that $S=P_\infty ^{(K)}$, $W\cong \mathbb{R}^{m_1}\times
\mathbb{C}^{m_2}$, and that
    $$
W^{(0)}=\bigoplus_{\substack{j=m_1+1,\dots ,m_1+m_2
    \\
|\theta _j|=1}}\mathbb{C}
    $$
is the subspace of $W\cong \mathbb{R}^m$ on which $A$ acts
isometrically. Since $\alpha $ is ergodic,
$\textup{dim}_\mathbb{R}(W^{(0)})\le
\textup{dim}_\mathbb{R}(W)-2=m-2$.

Take, for example, the irreducible ergodic and nonexpansive
automorphism $\alpha $ of $X=\mathbb{T}^4$ determined by the matrix
    $$
A=\left(
    \begin{smallmatrix}
\hphantom{-}0&1&0&0
    \\
\hphantom{-}0&0&1&0
    \\
\hphantom{-}0&0&0&1
    \\
-1&1&1&1
    \end{smallmatrix}
\right)\in\textup{GL}(4,\mathbb{Z}).
    $$
If $\theta >1$ is the dominant eigenvalue of $A$, then the algebraic
number field $K=\mathbb{Q}[\theta ]$ has two real places $v_1,v_2$
(corresponding to the real roots $\theta _1=\theta $ and $\theta
_2=\theta ^{-1}$ of the characteristic polynomial $f=u^4-u^3-u^2-u+1$
of $A$) and one complex place $v_3$ (corresponding to the two complex
roots $\theta _3$ and $\bar{\theta }_3$ of $f$ of absolute value
$1$). Then $S^{(0)}=\{v_3\}$, $W^{(0)}\cong K_{v_3}=\mathbb{C}$, and
the central subgroup $X^{(0)}\subset X$ of $\alpha $ is a densely
embedded copy of $\mathbb{C}$.

For another example of this form we take the automorphism $\alpha $
of $X=\mathbb{T}^6$ defined by the matrix
    $$
B=\left(
    \begin{smallmatrix}
\hphantom{-}0&1&0&0&0&0
    \\
\hphantom{-}0&0&1&0&0&0
    \\
\hphantom{-}0&0&0&1&0&0
    \\
\hphantom{-}0&0&0&0&1&0
    \\
\hphantom{-}0&0&0&0&0&1
    \\
-1&1&1&1&1&1
    \end{smallmatrix}
\right)\in\textup{GL}(6,\mathbb{Z})
    $$
with dominant eigenvalue $\theta >1$. The algebraic number field
$K=\mathbb{Q}[\theta ]$ has two real places $v_1,v_2$ (corresponding
to the real roots $\theta _1=\theta $ and $\theta _2=\theta ^{-1}$ of
the characteristic polynomial $f=u^6-u^5-u^4-u^3-u^2-u+1$ of $B$) and
two complex places $v_3,v_4$ (corresponding to the four complex roots
$\theta _3,\theta _4$ and $\bar{\theta }_3,\bar{\theta }_4$ of $f$ of
absolute value $1$). Then $S^{(0)}=\{v_3,v_4\}$, $W^{(0)}\cong
K_{v_3}\oplus K_{v_4}=\mathbb{C}^2$, and the central subgroup
$X^{(0)}\subset X$ of $\alpha $ is a densely embedded copy of
$\mathbb{C}^2$.

\smallskip (2) Let $f=5u^2 - 6u + 5\in R_1$, and let $\alpha =\alpha
_{R_1/(f)}$ be the automorphism of the compact connected abelian
group $X=X_{R_1/(f)}$ defined in \eqref{eq:ideal}. Since $f$ is
irreducible and all roots of $f$ have absolute value $1$ (they are of
the form $\theta =\frac35\pm i\cdot \frac45$), $\alpha $ is ergodic
and nonexpansive by Theorem \ref{t:irreducible}. If $\theta $ is a
root of $f$ and $K=\mathbb{Q}(\theta )$, then $P\subset P_f^{(K)}$,
$S^{(0)}=P_\infty ^{(K)}$, $W=W^{(0)}\times \prod_{v\in P}K_v$, where
$W^{(0)}\cong\mathbb{C}$ and $\prod_{v\in P}K_v$ is zero-dimensional,
and $X\cong W/L$ for some discrete co-compact $\bar{\beta
}$-invariant subgroup $L\subset W$. In this example the central
subgroup $X^{(0)}\subset X$ is a densely embedded copy of
$\mathbb{C}$.

\smallskip (3) Let $f=6u^4 + 3u^3 + 10u^2 + 6u + 6\in R_1$, and let
$\alpha =\alpha _{R_1/(f)}$ be the automorphism of the compact
connected abelian group $X=X_{R_1/(f)}$ defined in \eqref{eq:ideal}.
Again $f$ is irreducible, all roots of $f$ have absolute value $1$,
and $\alpha $ is ergodic and nonexpansive by Theorem
\ref{t:irreducible}. If $\theta $ is a root of $f$ and
$K=\mathbb{Q}(\theta )$, then $P\subset P_f^{(K)}$, $S^{(0)}=P_\infty
^{(K)}$, $W=W^{(0)}\times \prod_{v\in P}K_v$, where
$W^{(0)}\cong\mathbb{C}^2$ and $\prod_{v\in P}K_v$ is
zero-dimensional, and $X\cong W/L$ for some discrete co-compact
$\bar{\beta }$-invariant subgroup $L\subset W$. Here the central
subgroup $X^{(0)}\subset X$ is a densely embedded copy of
$\mathbb{C}^2$.
    \end{exas}

The group $W^{(0)}\cong \mathbb{C}^{|S^{(0)}|}$ in \eqref{eq:Sc0} is
an algebra with respect to coordinate-wise addition and
multiplication. We define a map $\iota _0\colon K\longrightarrow
W^{(0)}$ by setting
    $$
\iota (a)_v=
    \begin{cases}
\iota _v(a)&\textup{if}\enspace v\in S^{(0)},
    \\
0&\textup{if}\enspace v\in S\smallsetminus S^{(0)}
    \end{cases}
    $$
for every $a\in K$ (cf. \eqref{eq:iotav}), set
    \begin{equation}
    \label{eq:xi}
\xi =\iota _0(\theta ),
    \end{equation}
where $\theta $ is the algebraic number appearing in Theorem
\ref{t:irreducible2} and \eqref{eq:notation}, and denote by
    \begin{equation}
    \label{eq:Gamma}
\Gamma =\overline{\{\xi ^m:m\in\mathbb{Z}\}}
    \end{equation}
the closure of the multiplicative subgroup $\{\xi
^m:m\in\mathbb{Z}\}\subset W^{(0)}$. Then $\Gamma $ is a compact
abelian multiplicative subgroup of $W^{(0)}$. For every $\gamma
=(\gamma _v)\in \Gamma $ we denote by $M_\gamma \colon
W^{(0)}\longrightarrow W^{(0)}$ multiplication by $\gamma $, i.e.
    \begin{equation}
    \label{eq:Mgamma}
M_\gamma w=(\gamma _vw_v)
    \end{equation}
for every $w=(w_v)\in W^{(0)}$.

    \begin{prop}
    \label{prop:natural}
If $\alpha $ is totally irreducible, then for any two distinct
elements $v,v'\in S^{(0)}$, the natural projection of $\Gamma $ to
$K_v\oplus K_{v'}$ is surjective.
    \end{prop}

    \begin{proof}
Let $\xi_{vv'}$ be defined by
    $$
(\xi_{vv'})_{\nu} =\smash[t]{
    \begin{cases}
\theta &\text{for $\nu = v,v'$,}
    \\
0&\text{otherwise.}
    \end{cases}
}
    $$
Clearly the projection of $\Gamma$ to $K_v\oplus K_{v'}$ is equal to $\overline{\{\xi_{vv'}
^m:m\in\mathbb{Z}\}}$. Let $\xi _v = \iota _v(\theta )\in\mathbb{C}$
and $\xi _{v'} = \iota _{v'}(\theta )\in\mathbb{C}$. Since $v,v'\in
S^{(0)}$ we know that $|\xi _v|=|\xi _{v'}|=1$ (cf. \eqref{eq:iotav}).

In order to prove our claim it suffices to show that, for any nonzero
element $(m,m')\in \mathbb{Z}^2$,
    \begin{equation}
    \label{eq:multind}
\xi_v^m \xi_{v'}^{m'}\neq 1.
    \end{equation}
That $\xi_v^m \neq 1$ for $m\neq 0$ follows from ergodicity ($\xi _v$
is a root of an irreducible polynomial with integer coefficients,
which is  noncylotomic if $\alpha$ is to be ergodic). To prove \eqref
{eq:multind} for the case where both $m,m'\neq 0$, we note that since
$\xi_v$ and $\xi_{v'}$ are conjugate under the Galois group of the
splitting field of the polynomial $f$, we also have that $\xi_{v'}^m
= \xi_3^{m'}$ for some $\xi _3 \in \mathbb{C}$ with $f(\xi_3)=0$ (it
could be that $\xi _3 =\xi _v$). We can now apply the same argument
for $\xi_3$ and obtain that $\xi_3 ^ m = \xi_4 ^ {m'}$ for some root
$\xi_4\in\mathbb {C}$ of $f$, etc. Since $f$ has finitely many roots,
we will eventually get an equation of the form $\xi_j^{m^{k}} =\xi_j
^{(-m')^{k}}$ for some positive integer $k$ and some root $\xi _j$ of
$f$. As all roots of $f$ are conjugate under the Galois group, this
shows that
    $$
\theta ^{m^{k}} =  \theta ^{(-m')^{k}}.
    $$
If $m^k\neq (-m')^k$ then $\theta $ is a root of unity, which is a
contradiction to ergodicity. Otherwise $m =\pm m'$, and either
$\xi_v^m = \xi_{v'}^{m}$ or $\xi_v ^ m = \xi_{v'} ^ {-m}$.

First suppose that $\xi _v^m=\xi _{v'}^m$. Since $v $ and $v'$ are
inequivalent valuations, $\xi_v\neq \xi_{v'}$, and hence
$\xi_v\xi_{v'}^{-1}$ is a nontrivial root of unity, contrary to the
hypothesis that $\alpha $ is totally irreducible (cf. Theorem
\ref{t:irreducible}).

If $\xi _v^m=\xi _{v'}^{-m}$, then the complex conjugate $\xi
'=\overline{\xi _{v'}}$ of $\xi _{v'}$ is again a root of $f$
satisfying that $\xi _v^m={\xi '}^m$, and the same argument as above
shows that $\xi_v{\xi '}^{-1}$ is a nontrivial root of unity. Again
this violates the total irreducibility of $\alpha $.
    \end{proof}

\section{Conditional measures on the leaves of the central foliation}

We assume that $\alpha $ and $X$ are of the form \eqref{eq:notation}
and use the notation of \eqref{eq:Pc}--\eqref{eq:YL}. Write
$\mathcal{F}$ for the foliation of $X$ by the cosets of the central
subgroup $X^{(0)}=\pi (W^{(0)})\subset X$ (cf. \eqref{eq:Sc0}), and
fix a nonatomic $\alpha $-invariant probability measure $\mu $ on the
Borel field $\mathscr{S}=\mathscr{B}_X$ of $X$. Note that we do not
make any assumptions regarding ergodicity of $\mu$.

Since the central subgroup $X^{(0)}$ is dense by irreducibility, the
foliation of $X$ into cosets of $X^{(0)}$ has no Borel cross-section,
and one cannot generally decompose $\mu $ directly into a family of
measures supported on the individual leaves of $\mathcal{F}$. In
order to overcome this difficulty we break up each of these leaves
into countably many atoms of an appropriate sub-sigma-algebra
$\mathscr{A}\subset \mathscr{S}$, decompose the measure $\mu $ with
respect to this sigma-algebra, and re-combine the conditional
measures supported by the individual atoms on each leaf into a
leaf-measure.

It will be necessary to work not just with one such sigma-algebra
$\mathscr{A}$ but with a sequence $(\mathscr{A}^{(k)},\,k\ge 1)$ of
sigma-algebras whose atoms consist of larger and larger pieces of
leaves of $\mathcal{F}$. In order to describe these sigma-algebras we
fix an integer $q>1$ with $|q|_v=1$ for every $v\in P$ and set
$\Lambda =\frac1qL\subset K$. Then $\iota (\Lambda )$ is a discrete
co-compact subgroup of $W$ (cf. \eqref{eq:ic}), and we choose a Borel
set $\Delta _0\subset W$ with compact closure such that
    \begin{equation}
    \label{eq:Delta}
    \begin{gathered}
\Delta \cap (\Delta +\iota (a))=\varnothing \enspace \textup{for
every nonzero}\enspace a\in\Lambda ,
    \\
\bigcup_{a\in\Lambda }\Delta +\iota (a)=W.
    \end{gathered}
    \end{equation}
The first equation in \eqref{eq:Delta} implies that the restriction
to $\Delta $ of the map $\pi \colon W\longrightarrow X$ in
\eqref{eq:pi} is injective, and that $\pi (\Delta )$ is therefore a
Borel subset of $X$. After replacing $\Delta $ by $\Delta +w_0$ for
some $w_0\in W$, if necessary, we may take it that the sets
    \begin{equation}
    \label{eq:Q}
\mathcal{Q}=\{\pi (\Delta +\iota (a)):a\in\Lambda \}
    \end{equation}
form a Borel partition of $X$ into $N=|L/qL|$ sets with the following
properties.
    \begin{enumerate}
    \label{enum}
    \item[(i)]
$\mu (\partial Q)=0$ for every $Q\in\mathcal{Q}$;
    \item[(ii)]
For every $a\in \Lambda $, the restriction of the map $\pi $ in
\eqref{eq:pi} to $\Delta  + \iota (a)$ is injective, and $\pi (\Delta
+ \iota (a))=\pi (\Delta  + \iota (a'))$ if and only if $a-a'\in L$;
    \item[(iii)]
For every $Q\in\mathcal{Q}$ and $w\in W$, the set $W^{(0)}\cap(\pi
^{-1}(Q)-w)\subset W^{(0)}\cap \bigcup_{a\in L}(\mathcal{C}-w)$ is a
countable union of sets with disjoint and compact closures.
    \end{enumerate}

Let $T_y$ denote the map
    \begin{equation}
    \label{eq:T}
T_yx=x+y
    \end{equation}
for every $x,y\in X$. We denote by $B_{W^{(0)}}(w,r)$ the ball of
radius $r>0$ around $w$ in $W^{(0)}$; while it will not be important
for us which norm we use in $W^{(0)}$, the natural norm to take is
    \begin{equation}
    \label{eq:norm2}
\|w\| =\max_{v\in S^{(0)} } |w_v|
    \end{equation}
(cf. \eqref{eq:norm}). Finally we write $B_{\mathcal{F}}(x,r)$ for
the ball of radius $r$ around $x$ in the leaf $x+\pi(W^{(0)})$ of
$\mathcal{F}$, i.e.
    \begin{equation}
    \label{eq:BF}
B_{\mathcal{F}}(x,r)=x+\pi(B_{W^{(0)}}(0,r))=T_x\circ \pi(B_{W^{(0)}}(0,r)).
    \end{equation}

    \begin{prop}
    \label{p:invariance}
There exist a sequence of fundamental domains
$(\Delta^{(n)},\,n\linebreak[0]\ge1)$ for $\iota (\Lambda)$ and a
corresponding sequence of partitions $(\mathcal{Q}^{(n)},\,n\ge1)$ of
$X$ in \eqref{eq:Q} with the properties \textup{(i)--(iii)}
\vpageref[above]{enum}, such that
    \begin{equation}
    \label{eq:invariance}
\sum_{Q\in\mathcal{Q}^{(n)}}\mu \bigl((Q +
\pi(B_{W^{(0)}}(0,n)))\bigtriangleup Q\bigr)\le 2^{-n}
    \end{equation}
for every $n\ge1$.
    \end{prop}

    \begin{proof}
By choosing, for every $n\ge1$, an appropriate fundamental domain
$\Delta _n\subset W$ for $\Lambda $ with the properties
\eqref{eq:Delta} we can construct a sequence
$(\mathcal{Q}_n,\,n\ge1)$ of Borel partitions \eqref{eq:Q} satisfying
the conditions (i)--(iii) \vpageref[above]{enum} such that
    $$
\sum_{Q\in\mathcal{Q}_n}\lambda_X \bigl((Q +
\pi(B_{W^{(0)}}(0,n)))\bigtriangleup Q\bigr)\le 2^{-n}
    $$
for every $n\ge1$. Since for any two Borel sets $Q,B\subset X$,
    \begin{align*}
\int\mu \bigl((Q+s+B)&\bigtriangleup (Q+s)\bigr)\,d\lambda_X(s)
    \\
&=\smash[t]{\iint \bigl|1_{Q+s+B}(x)-1_{Q+s}(x)\bigr|\, d\mu
(x)\,d\lambda_X(s)}
    \\
&=\lambda_X\bigl((Q+B)\bigtriangleup Q\bigr),
    \end{align*}
where $1_{Q+s+B}$ and $1_{Q+s}$ are the indicator function of the
sets $Q+s+B$ and $Q+s$, there is a sequence $(x_n,\,n\ge1)$ so that
the translated partitions $\mathcal{Q}^{(n)} = \mathcal{Q}_n +x_n$
satisfy \eqref{eq:invariance} and the conditions \textup{(i)--(iii)}
for every $n\ge1$. For later use we choose a bounded sequence
$(w_n,\,n\ge1)$ in $W$ with $\pi (w_n)=x_n$ for every $n\ge1$ and set
$\Delta ^{(n)}=\Delta _n+w_n,\,n\ge1$.
    \end{proof}

    \begin{defi}
    \label{d:equivalent}
Let $\mathscr{A}\subset \mathscr{S}$ be a countably generated
sigma-algebra, and let $\mathcal{C}\subset \mathscr{A}$ be a
countable algebra which generates $\mathscr{A}$. The \emph{atom}
$[x]_\mathscr{A}$ of a point $x\in X$ in $\mathscr{A}$ is defined as
    $$
[x]_\mathscr{A}=\bigcap_{C\in\mathcal{C}:x\in
C}C=\bigcap_{A\in\mathscr{A}:x\in A}A.
    $$
    \end{defi}

    \begin{lemm}
    \label{l:algebra}
For every $n\ge1$, let $\Delta ^{(n)}$ be the fundamental domain for
$\iota (\Lambda )\subset W$ and $\mathcal{Q}^{(n)}$ the partition of
$X$ described in Proposition \ref{p:invariance}. Then there exist a
countably generated sigma-algebra $\mathscr{A}^{(n)}\subset
\mathscr{S}$ with
    \begin{equation}
    \label{eq:Ax}
[x]_{\mathscr{A}^{(n)}}=\pi \bigl((\Delta ^{(n)}+\iota
(a))\cap(W^{(0)}+w)\bigr)
    \end{equation}
for every $x\in X$, where $[x]_{\mathscr{A}_{\mathcal{Q}^{(n)}}}$ is
the atom of $\mathscr{A}^{(n)}$ containing $x$ and $a\in \Lambda $
and $w\in W$ satisfy that $\pi (w)=x$ and $w\in\Delta ^{(n)}+\iota
(a)$.
    \end{lemm}

    \begin{proof}
Fix $a\in \Lambda $ for the moment. We set $W'=\bigl(\prod_{v\in
S\smallsetminus S^{(0)}}K_v\bigr)$, denote by $\kappa \colon
W=W^{(0)}\times W'\longrightarrow W'$ the second coordinate
projection (cf. \eqref{eq:Sc0}), and write $\mathscr{B}_{W'}$ for the
(countably generated) Borel field of $W'$. The sigma-algebra
$\mathscr{A}=\{\kappa ^{-1}(B)\cap (\Delta +\iota (a)):B\in
\mathscr{B}_{W'}\}$ of subsets of $\Delta +\iota (a)$ is again
countably generated, and its atoms are of the form $(\Delta +\iota
(a))\cap (W^{(0)}+w),\,w\in \Delta +\iota (a)$. Since the restriction
of $\pi $ to $\Delta +\iota (a)$ is injective, $\pi $ maps
$\mathscr{A}$ to a countably generated sigma-algebra $\mathscr{A}_Q$
of subsets of $Q=\pi (\Delta +\iota (a))\in\mathcal{Q}^{(n)}$ whose
atoms are of the required form. The sigma-algebra $\mathscr{A}^{(n)}$
is defined as the unique sub-sigma-algebra of
$\mathscr{S}=\mathscr{B}_X$ which contains the partition
$\mathcal{Q}^{(n)}$ and induces $\mathscr{A}_Q$ on each
$Q\in\mathcal{Q}^{(n)}$.
    \end{proof}

For any countably generated sigma-algebra $\mathscr{A}
\subset\mathscr{S}$ we consider the decomposition of $\mu $ with
respect to the sigma-algebra $\mathscr{A}$, i.e. a set of probability
measures $\{\mu _x^{\mathscr{A}}: x\in X\}$ on $X$ with the following
properties.
    \begin{enumerate}
    \item
    \label{cond}
For all $x,x'\in X$ with $[x]_{\mathscr{A}}= [x']_{\mathscr{A}}$,
    $$
\mu _x^{\mathscr{A}}=\mu _{x'}^{\mathscr{A}}\enspace
\textup{and}\enspace \mu _x^{\mathscr{A}} ([x]_{\mathscr{A}})=1,
    $$

    \item
For every $B\in\mathscr{S}$, the map $x\mapsto \mu
_x^{\mathscr{A}}(B)$ is Borel (and hence $\mathscr{A}$-measurable),
    \item
For every bounded Borel map $f\colon X\longrightarrow \mathbb{R}$,
    $$
\int f\,d\mu _x^\mathcal{A}=E_\mu (f|\mathcal{A})(x)
    $$
for $\mu\textsl{-a.e.}\;x\in X$, where $E_\mu(\cdot |\cdot )$ denotes
conditional expectation.
    \end{enumerate}

In order to make notation less cumbersome we set, for every
$n\in\mathbb{Z}$ and $k\ge 1$,
    \begin{equation}
    \label{eq:awkward}
\mathscr{A}_n^{(k)}=\alpha ^{-n}(\mathscr{A}_{\mathcal{Q}^{(k)}}),
\qquad \mathscr{A}^{(k)}=\mathscr{A}_0^{(k)}
=\mathscr{A}_{\mathcal{Q}^{(k)}},
    \end{equation}
and denote by $\{\mu _x^{\mathscr{A}^{(k)}_n}:x\in X\}$ and $\{\mu
_x^{\mathscr{A}^{(k)}\vee \mathscr{A}_n^{(k)}}:x\in X\}$ the
decompositions of $\mu $ with respect to the sigma-algebras
$\mathscr{A}_n^{(k)}$ and $\mathscr{A}^{(k)}\vee
\mathscr{A}^{(k)}_n$, respectively.

    \begin{defi}
    \label{d:localfinite}
A Borel measure $\rho $ on $W^{(0)}$ is \textsl{locally finite} if
$\rho (C)<\infty $ for every compact set $C\subset W^{(0)}$. Let
$M_\infty (W^{(0)})$ be the set of all locally finite (and hence
sigma-finite) Borel measures on $W^{(0)}$, furnished with the
smallest topology in which the map $\rho \mapsto \int f\,d\rho $ from
$M_\infty (W^{(0)})$ to $\mathbb{R}$ is continuous for every
continuous map $f\colon W^{(0)}\longrightarrow \mathbb{R}$ with
compact support. In this topology $M_\infty (W^{(0)})$ is a separable
metrizable space.
    \end{defi}

For every $a\in W^{(0)}$ we denote by
    \begin{equation}
    \label{eq:Ta}
\bar{T}_wv=v+w,\;w\in W^{(0)},
    \end{equation}
the translation by $w$ on $W^{(0)}$. The maps $\rho \mapsto \rho
\bar{T}_w$ and $\rho \mapsto \rho \bar{\beta }$ are homeomorphisms of
$M_\infty (W^{(0)})$ for every $w\in W^{(0)}$, where $\bar{\beta }$
is defined in \eqref{eq:barbetac}.

For the next theorem, we take $r_0$ to be large enough so that
    \begin{equation}
    \label{eq:r_0 definition}
[x]_{\mathscr{A}^{(1)}}\subset B_{\mathcal{F}}(x,r_0)
    \end{equation}
for all $x\in X$ (cf. \eqref{eq:BF}).

    \begin{prop}
    \label{p:rho}
There is a Borel map $x\mapsto \rho _x$ from $X$ to $M_\infty
(W^{(0)})$ and an $\alpha $-invariant Borel set $X'$ of full $\mu
$-measure with the following properties \textup{(}for notation we
refer to \eqref{eq:barbetac}, \eqref{eq:awkward} and
\eqref{eq:Ta}\textup{)}.
    \begin{enumerate}
    \item
For every $x\in X'$, every bounded Borel set $B\subset W^{(0)}$ and
every sufficiently large $k$,
    \begin{equation}
    \label{eq:rho}
\rho _x (B)= \frac{1}{\mu _x^{\mathscr{A}^{(k)}}
(B_{\mathcal{F}}(x,r_0))}\mu _x^{\mathscr{A}^{(k)}} (T_x\circ \pi
(B));
    \end{equation}
    \item
For every $x\in X$,
    \begin{equation}
    \label{eq:rhoalphabeta}
\rho _x=\rho _{\alpha x}\bar{\beta };
    \end{equation}
    \item
There exists a Borel map $K_\mu \colon X\times W^{(0)}\longrightarrow
\mathbb{R}$ so that, for every $x\in X'$ and every $w\in W^{(0)}$
with $x+\pi (w)\in X'$,
    \begin{equation}
    \label{eq:mutranslates}
e^{K_\mu (x,w)} \rho _{x-\pi (w)}=\rho _{x}\bar{T}_w.
    \end{equation}
    \end{enumerate}
    \end{prop}

We begin the proof of Proposition \ref{p:rho} with a lemma.

    \begin{lemm}
    \label{l:bigballs}
There exists a Borel set $X'\subset X$ with $\mu (X')=1$, which is
invariant under $\bar T_w$ for every $w\in W^{(0)}$, so that for all
$x\in X'$ and $r>0$,
    \begin{equation}
    \label{eq:bigballs}
B_{\mathcal{F}}(x,r) \subset [x]_{\mathscr{A}^{(k)}}
    \end{equation}
for every sufficiently large $k$.
    \end{lemm}

    \begin{proof}
The set
    $$
N_{r,k}=\{x\in X: B_{\mathcal{F}}(x,r) \not \subset [x]_{\mathscr{A}^{(k)}}\}.
    $$
is equal to
    \begin{align*}
\smash{\bigcup_{\substack{Q,Q'\in\mathcal{Q}^{(k)}
    \\
Q\ne Q'}}} \bigl(Q\cap(Q'+\pi(B_{W^{(0)}} (0,r)))\bigr)&=
\bigcup_{Q'\in\mathcal{Q}^{(k)}} \bigl((Q'+\pi(B_{W^{(0)}} (0,r)))
\smallsetminus Q'\bigr)
    \\
&\subset \bigcup_{Q'\in\mathcal{Q}^{(k)}} \bigl((Q'+\pi(B_{W^{(0)}}
(0,r)))\bigtriangleup Q'\bigr).
    \end{align*}
and hence Borel. Since $\sum_{k\ge1} \mu (N_{r,k}) <\infty$ for every
$r>0$ by \eqref{eq:invariance}, it follows that
    $$
X'=X\smallsetminus \bigcap_{r>0} \bigcup_{n\ge 1}\bigcap_{k\geq n} N_{r,k}
    $$
is a Borel set of full measure. From the definition of $X'$ it is
also clear that any $x\in X'$ satisfies \eqref{eq:bigballs}, and that
$X'$ consists of a union of full $\mathcal{F}$ leaves, i.e. that it
is invariant under $\bar T_w$ for any $w\in W^{(0)}$.
    \end{proof}

    \begin{proof}[Proof of Proposition \ref{p:rho}]
We take $X'$ to be the set of all $x\in X$ with the following properties.
    \begin{enumerate}
    \label{enum3}
    \item
For every $r>0$ and $n\in\mathbb{Z}$, and for every sufficiently
large $k\ge 1$ (depending on $r$ and $n$),
    \begin{equation}
    \label{eq:bigballstwo}
B_{\mathcal{F}}(x,r) \subset[x]_{\mathscr{A}^{(k)}_n};
    \end{equation}
    \item
For every $k,l\ge1$ and $n\in\mathbb{Z}$,
    \begin{gather}
    \label{eq:l,kone}
\mu _x^{\mathscr{A}^{(k)}}([x]_{\mathscr{A}^{(k)} \vee
\mathscr{A}^{(l)}_n})>0, \qquad \mu _x^{\mathscr{A}^{(l)}_n}
([x]_{\mathscr{A}^{(k)} \vee \mathscr{A}^{(l)}_n})>0,
    \\
    \label{eq:l,ktwo}
    \begin{aligned}
\mu _x^{\mathscr{A}^{(k)} \vee \mathscr{A}^{(l)}_n}&=\frac1{\mu
_x^{\mathscr{A}^{(k)}}([x]_{\mathscr{A}^{(k)} \vee
\mathscr{A}^{(l)}_n})} \cdot \mu
_x^{\mathscr{A}^{(k)}}\hspace{-1.5mm} \bigm|_{[x]_{\mathscr{A}^{(k)}
\vee \mathscr{A}^{(l)}_n}}
    \\
&= \frac1{\mu _x^{\mathscr{A}^{(l)}_n}([x]_{\mathscr{A}^{(k)} \vee
\mathscr{A}^{(l)}_n})} \cdot \mu _x^{\mathscr{A}^{(l)}_n}
\hspace{-1.5mm}\bigm|_{[x]_{\mathscr{A}^{(k)} \vee
\mathscr{A}^{(l)}_n}},
    \end{aligned}
    \end{gather}
where $\mu _x^{\mathscr{A}^{(k)}} \hspace{-1.5mm}\bigm|_{[x]_{\mathscr{A}^{(k)}
\vee \mathscr{A}^{(l)}_n}}$ and $\mu _x^{\mathscr{A}^{(l)}_n}
\hspace{-1.5mm}\bigm |_{[x]_{\mathscr{A}^{(k)} \vee
\mathscr{A}^{(l)}_n}}$ are the restrictions of $\mu
_x^{\mathscr{A}^{(k)}}$ and $\mu _x^{\mathscr{A}^{(l)}_n}$ to the
atom $[x]_{\mathscr{A}^{(k)} \vee \mathscr{A}^{(l)}_n}$ of $x$ in
$\mathscr{A}^{(k)} \vee \mathscr{A}^{(l)}_n$;
    \item
For every $k\ge1$ and $n\in\mathbb{Z}$,
    \begin{equation}
    \label{eq:transform}
\mu _x^{\mathscr{A}^{(k)}_n} =\mu _{\alpha ^nx}^{\mathscr{A}^{(k)}} \alpha ^n.
    \end{equation}
    \end{enumerate}
Note that by \eqref{eq:r_0 definition} and \eqref{eq:l,kone} (with
$l=1$ and $n=0$),
    $$
\mu _x^{\mathscr{A}^{(k)}}(B_{\mathcal{F}}(x,r_0))>0
    $$
for every $k\ge1$ and $x\in X'$. Furthermore, by \eqref
{eq:bigballstwo} and \eqref {eq:l,ktwo} (again with $n=0$),
    $$
\frac{1}{\mu _x^{\mathscr{A}^{(k)}} (B_{\mathcal{F}}(x,r_0)) }\mu
_x^{\mathscr{A}^{(k)}} \hspace{-1.5mm}\bigm|_{[x]_{\mathscr{A}^{(k)}
\vee \mathscr{A}^{(l)}}}= \frac{1}{\mu _x^{\mathscr{A}^{(l)}}
(B_{\mathcal{F}}(x,r_0))}\mu _x^{\mathscr{A}^{(l)}}
\hspace{-1.5mm}\bigm|_{[x]_{\mathscr{A}^{(k)} \vee \mathscr{A}^{(l)}}}
    $$
for all $x\in X'$ and all sufficiently large $k,l$, so that
    \begin{equation}
    \label{eq:rhox}
\rho _x (B)=\lim_{k\to\infty} \frac{1}{\mu _x^{\mathscr{A}^{(k)}}
(B_{\mathcal{F}}(x,r_0))}\mu _x^{\mathscr{A}^{(k)}}(T_x \circ \pi (B))
    \end{equation}
exists for every $x\in X'$ and every Borel set $B\subset W^{(0)}$.

Equation \eqref{eq:mutranslates} easily follows from the fact that,
for every $x\in X'$, every $w\in W^{(0)}$ with $x-\pi (w)\in X'$, and
every sufficiently large $k$,
    $$
y=x-\pi(w)\in [x]_{\mathscr{A}^{(k)}},
    $$
and hence
    $$
\mu _y^{\mathscr{A}^{(k)}}= \mu _x^{\mathscr{A}^{(k)}}.
    $$
The sequence
    $$
\log\,\frac{\mu _{x-\pi (w)}^{\mathscr{A}^{(k)}}
(B_{\mathcal{F}}(x-\pi (w),r_0))}{\mu _x^{\mathscr{A}^{(k)}}
(B_{\mathcal{F}}(x,r_0))}= \log\,\frac{\mu _{x}^{\mathscr{A}^{(k)}}
(B_{\mathcal{F}}(x-\pi (w),r_0))}{\mu _x^{\mathscr{A}^{(k)}}
(B_{\mathcal{F}}(x,r_0))}
    $$
is eventually constant, and we set
    $$
K_\mu (x,w)=
    \begin{cases}
\smash[b]{\lim\limits_{k\to\infty }\,\log\,\frac{\mu
_{x}^{\mathscr{A}^{(k)}} (B_{\mathcal{F}}(x-\pi (w),r_0))}{\mu
_x^{\mathscr{A}^{(k)}}
(B_{\mathcal{F}}(x,r_0))}}&\textup{for}\enspace x\in X'\enspace
\textup{and}\enspace w\in W^{(0)}
    \\
&\qquad \textup{with}\enspace x-\pi (w)\in X',
    \\
0&\textup{otherwise}.
    \end{cases}
    $$
Equation \eqref {eq:rhoalphabeta} is immediate from \eqref
{eq:bigballstwo}, \eqref {eq:l,ktwo} and \eqref{eq:transform} (with
$k=l$ and $n=1$).

Finally we extend the map $x\mapsto \rho _x$ to $X$ by setting $\rho
_x=0$ for every $x\in X\smallsetminus X'$ and note that the resulting
map from $X$ to $M_\infty (W^{(0)})$ is Borel.
    \end{proof}

\section{Finiteness of the central leaf measures}
    \label{s:finite}

    \begin{theo}
    \label{t:finite}
Let $\alpha$ be a nonexpansive, ergodic and totally irreducible
automorphism of a compact connected abelian group $X$ with normalized
Haar measure $\lambda _X$, and let $\mu $ be an $\alpha $-invariant
probability measure on $X$ which is singular with respect to $\lambda
_X$. Then there exists a Borel set $X'\subset X$ with $\mu (X')=1$
such that $\rho _x(W^{(0)})<\infty$ for every $x\in X'$ \textup{(}cf.
\eqref{eq:rho}\textup{)}.
    \end{theo}

We begin the proof of Theorem \ref{t:finite} with a series of lemmas
in which we denote the $l$-th derivative of a map $f$ by $f^{(l)}$.

    \begin{lemm}
    \label{l:poly}
For every $s\ge1$ we can find a constant $A_s>0$ such that, for every
polynomial $p(x) = \sum_{l=0}^{2s-1} a_lx^l$ of degree $\leq 2s-1$
and every $\varepsilon > 0$,
    $$
\sup_{t\in(-\varepsilon ,\varepsilon )}|p(t)| \geq A_s \cdot
\max_{0\le l\leq 2s-1} (\varepsilon ^l|a_l|).
    $$
    \end{lemm}

    \begin{proof}
The statement of the lemma is clearly unchanged by rescaling $p$ and
$\varepsilon $, so that we may assume that $\varepsilon = 1$ and
$\max_l |a_l| = 1$. We can now set
    \begin{equation}
A_s = \inf \biggl\{\sup_{t\in(-1,1)}|p(t)|: p(x)=\sum_{l=0}^{2s-1}
a_lx^l \enspace \textup{with}\enspace \max_l |a_l| = 1 \biggr\} >
0.\tag*{$\qed$}
    \end{equation}
\renewcommand{\qed}{}
    \end{proof}

    \begin{lemm}
    \label{l:sublemma}
Let $\varepsilon >0$, $s\ge1$, and let $A_s>0$ be the constant
appearing in Lemma \ref{l:poly}. Then
    \begin{equation}
    \label{eq:sublemma}
\sup_{t\in (-\varepsilon ,\varepsilon )}|f(t)|\ge \frac{A_sB}{2(2s-1)!}.
    \end{equation}
for every $B>0$ and every map $f\colon (-\varepsilon ,\varepsilon
)\longrightarrow \mathbb{R}$ with $2s$ derivatives at every point
such that
    \begin{equation}
    \label{eq:sublemma cond}
\smash[b]{\max_{0\le l\le 2s-1}|f^{(l)}(t)| \ge \frac{B}{\varepsilon
^l}}\enspace \enspace \textup{for every}\enspace t\in(-\varepsilon
,\varepsilon )
    \end{equation}
and
    \begin{equation}
    \label{eq:remainder}
\sup_{t\in(-\varepsilon ,\varepsilon )}|f^{(2s)}(t)|<
\frac{A_sB}{\varepsilon ^{2s}}.
    \end{equation}
    \end{lemm}

    \begin{proof}
Consider the Taylor expansion
    $$
p(x)=\sum_{l=0}^{2s-1} \frac{f^{(l)}(0)}{l!} x^l
    $$
of $f$ of degree $2s-1$. From Lemma \ref{l:poly} we know that there
is some $t\in(-\varepsilon ,\varepsilon )$ with $|p(t)|\geq
\frac{A_sB}{(2s-1)!}$, and Taylor's Theorem allows us to find a $\xi
\in [0,1]$ with
    $$
\smash[b]{f(t)=p(t) + \frac{f^{(2s)}(\xi t)}{(2s)!}t^{2s}.}
    $$
Thus
    \begin{align}
|f(t)|&\geq |p(t)| - \varepsilon ^{2s}\cdot \bigl(
\sup\nolimits_{t\in(-\varepsilon ,\varepsilon )} |f^{(2s)}(t)| \bigr)
\bigm/(2s)! \notag
    \\
&\ge  \frac{A_s B}{(2s-1)!} - \frac{A_sB}{(2s)!} \ge
\frac{A_sB}{2(2s-1)!}. \tag*{$\qed$}
    \end{align}
\renewcommand{\qed}{}
    \end{proof}

    \begin{lemm}
    \label{l:calculation}
Let $p(t)=\sum_{k=1}^s (a_k \cos(2\pi m_k t)+b_k \sin(2\pi m_k t))$
be a trigonometric polynomial, where the $m_k,\,k=1,\dots ,s$, are
distinct positive integers. Let $\|p\| = \max_{k=1,\dots ,s}
(|a_k+ib_k|)$. Then there exists a constant $c_2>0$, which depends on
$s$ and $M=\max_{k=1,\dots ,s} |m_k|$, but not on the coefficients
$a_k,b_k,\,k=1,\dots ,s$, such that
    \begin{equation}
    \label{eq:calculation}
\biggl|\int^1_0 e^{ip(t)}\,dt \biggr|\le c_2\cdot  \|p\|^{-1/2s}.
    \end{equation}
    \end{lemm}

    \begin{proof}
We first claim that, unless all coefficients $a_k$ and $b_k$ are $0$,
the derivative $p'$ of $p$ does not have zeros of order $> 2s-1$.
Indeed,
    $$
p^{(2l)}(t) = \sum_{k=1}^s(-1)^l (2\pi m_k)^{2l}(a_k\cos(2\pi m_k t)
+ b_k\sin(2\pi m_k t))
    $$
for every $l\ge 0$. If $p^{(2l)}(t_0) = 0$ for $l=1,\dots ,s$, the
nonsingularity of the Vandermonde matrix (due to our hypothesis that
the $m_k$ are all distinct) implies that
    $$
a_k\cos(2\pi m_k t_0)+b_k \sin(2\pi m_k t_0)  = 0
    $$
for $k = 1,\dots ,s$. Similarly, if $p^{(2l-1)}(t_0) = 0$ for $l =
1,\dots ,s$, then
    $$
-a_k \sin(2\pi m_k t_0) + b_k \cos(2\pi m_k t_0) = 0
    $$
for $k=1,\dots ,s$, and by combining these statements we get that
$a_k=b_k=0$ for $k=1,\dots ,s$. In fact, this argument gives more:
since one can bound the norm of the inverse of the Vandermonde matrix
for all choices $0<m_1<\dots<m_s\leq M$ by some function of $M$,
there exists a constant $c_M'>0$ depending only on $M$, such that
    $$
\max_{1\leq l\leq 2s} |p^{(l)}(t)|\ge c_M' \|p\|
    $$
for every $t\in\mathbb{R}$ and every choice of the coefficients
$a_k,b_k$ in $p$.

Trivially there exists, for every $l\ge 0$ and $M\ge1$, a constant
$c_{l,M}'>0$ such that
    $$
|p^{(l)}(t)|\le c_{l,M}'\|p\|
    $$
for every $l\ge0$ and $t\in\mathbb{R}$.

In order to complete the proof of Lemma \ref{l:calculation} we recall
the van der Corput Lemma in \cite[p. 220]{Katznelbook}: if $\phi$ is
a real-valued function on an interval $[a,b]\subset \mathbb{R}$ with
a monotonic derivative satisfying that $\phi'(t)>A>0$ for every $t\in
[a,b]$, then
    $$
\left|\int^b_a e^{i\phi(t)}\,dt\right| \leq \frac{4}{A}.
    $$
Since a trigonometric polynomial of degree $M$ such as $p''(t)$ can
have at most $2M$ roots in the interval $[0,1)$, the interval $[0,1]$
can be divided into at most $2M+1$ subintervals $I_1,I_2,\dots$, on
each of which $p'$ is monotonic. By applying the van der Corput Lemma
on each of these subintervals separately we have that, for any $A>0$,
    \begin{equation}
    \label{eq:any A}
\left|\int^1_0 e^{ip(t)}\,dt \right |\le \frac{8M+4}{A}+\lambda
(\{0\leq t\leq 1: |p'(t)|<A\}),
    \end{equation}
where $\lambda $ is the Lebesgue measure on $\mathbb{R}$.

It remains to estimate $\lambda (\{0\leq t\leq 1: |p'(t)|<A\})$. In
fact, we claim that there exists a constant $c''>0$ with
    \begin{equation}
    \label{eq:claim measure}
\lambda (\{0\leq t\leq 1: |p'(t)|<A\})\le
c''\left(\frac{A}{\|p\|}\right )^{\frac{1}{2s-1}}
    \end{equation}
for every $A>0$. Estimates of this kind can be found e.g. in
\cite{KM}; for completeness we provide a proof below.

As $p'$ is monotonic on every subinterval $I_k$,
    $$
I_k'=I_k\cap \{0\leq t\leq 1: |p'(t)|<A\}
    $$
is connected and hence an interval, and we apply Lemma
\ref{l:sublemma} with $f=p'$ on some sufficiently small subinterval
$I_k''\subset I'_k$. The conditions \eqref{eq:sublemma cond} and
\eqref{eq:remainder} are clearly satisfied for some $B=\|p\|\cdot
\lambda (I_k'')^{2s-1}$, and Lemma \ref{l:sublemma} guarantees the
existence of a constant $c_1>0$ with
    $$
\smash[b]{A\ge \sup_{t\in I'_k} p'(t)\ge c_1\|p\|\lambda (I'_k)^{2s-1}}
    $$
or
    $$
\smash[t]{\lambda (I'_k)\le c_1^{-\frac{1}{2s-1}}\cdot
\biggl(\frac{A}{\|p\|} \biggr) ^{\frac{1}{2s-1}}.}
    $$
By summing over $k$ we obtain \eqref{eq:claim measure}.

According to \eqref{eq:any A} and \eqref{eq:claim measure},
    $$
\left|\int^{2\pi}_0 e^{2\pi ip(t)}\,dt \right |\le
\frac{2}{A}+c''\cdot \biggl(\frac{A}{\|p\|} \biggr) ^{\frac{1}{2s-1}},
    $$
and by taking $A = \|p\|^{1/2s}$ we get \eqref{eq:calculation}.
    \end{proof}

We derive from this the following estimate.

    \begin{lemm}
    \label{l:character inequality}
For every nontrivial character $a\in\hat{X}$ there exists a constant
$c_a>0$ with
    $$
\int_{\Gamma} \langle a,\pi (M_{\gamma } w)\rangle \,d\gamma \le c_a
\cdot \min(1,\|w\|^{-1/2s})
    $$
for every $w\in W^{(0)}$, where $s=|S^{(0)}|$, and where $\Gamma $
and $M_\gamma $ are defined in \eqref{eq:Gamma} and \eqref{eq:Mgamma}.
    \end{lemm}

    \begin{proof}
We recall that $W=\prod_{v\in S} K_v$, consider each $K_v$ (by abuse
of notation) as an additive subgroup of $W$, and identify $X$ with
$W/\iota (L)$ as in \eqref{eq:YL}. Let $a$ be a nontrivial character
of $X=W/\iota (L)$, and let $f\colon w\mapsto f(w) = \langle a,\pi
(w)\rangle $ be the corresponding character of $W$. We write $f_0 =
f|_{W^{(0)}}$ for the restriction of $f$ to $W^{(0)}$. The
isomorphisms $\iota _v\colon K_v \longrightarrow \mathbb{C},\,v\in
S^{(0)}$, in \eqref{eq:iotav} allow us to write $f_0\colon
W^{(0)}\longrightarrow \mathbb{C}$ as the map
    $$
w\mapsto f_0 (w) = e^{2\pi i\sum_{v\in S^{(0)}} \Re(a_v \iota _v(w_v))},
    $$
where $a_v\in \mathbb{C}$ for every $v\in S^{(0)}$, and where $\Re$
denotes the real part. Since the image $\pi(K_v)$ of $K_v$ is dense
in $X$ for every $v$ by irreducibility, $f_0|_{K_v}$ is a nontrivial
character, hence $a_v\neq 0$ for every $v\in S^{(0)}$.

Let $v,v'$ be distinct elements of $S^{(0)}$. By Proposition
\ref{prop:natural}, the projection of $\Gamma $ to $K_v \oplus
K_{v'}\subset W^{(0)}$ maps $\Gamma $ onto the set
    $$
\{w \in W: |w_v| = |w_{v'}| =1\enspace \text{and all other
coordinates are 0}\}.
    $$
Hence there exists a closed one-dimensional subgroup
    $$
\Gamma _0=\bigl\{\bigl(\iota_v^{-1} (z^{m_v})\bigr)_{v\in S^{(0)}}:
|z|=1\bigr\}\subset \Gamma \subset W^{(0)}
    $$
such that the integers $m_v,\,v\in S^{(0)}$, are all distinct.

Now we use Lemma \ref{l:calculation} to check that there exists a
constant $c$ with
    $$
\biggl|\int_{\Gamma_0}f(M_{\gamma}w)\,d\gamma\biggr| =
\biggl|\int_{0}^1e^{2\pi i\Re\bigl(\sum_{v\in S^{(0)}}
a_v\iota_v(w_v)e^{2\pi im_vt}\bigr)}dt\biggr|\le c\cdot \|w\|^{-1/2s}
    $$
for every $w\in W^{(0)}$, and by integrating over $\Gamma $ we see that
    \begin{align}
\biggl|\int_{\Gamma}f(M_{\gamma} w)\,d\gamma\biggr| &=
\biggl|\int_{\Gamma} \int_{\Gamma_0} f(M_{\gamma_0} M_{\gamma}w)\,
d\gamma_0 \,d\gamma\biggr|\notag
    \\
&\le c\cdot \int_{\Gamma} \|M_{\gamma} w\|^{-1/2s} d\gamma = c\cdot
\|w\|^{-1/2s}.\tag*{$\qed$}
    \end{align}
\renewcommand{\qed}{}
    \end{proof}

    \begin{lemm}
    \label{l:character}
For every nontrivial character $a\in \hat{X}$ there exists a constant
$c_a>0$ such that
    \begin{equation}
    \label{eq:harmonic estimate}
    \begin{aligned}
\int _{\Gamma}\biggl| \int _{W^{(0)}} &\langle a,\pi (M_{\gamma
}w)\rangle \, d\bar{\tau } (w) \biggr|^2 d\gamma
    \\
&\le c_a\cdot \smash{\int \min (1,\|w-w'\|^{-1/2s})\, d\bar{\tau
}(w)\, d\bar{\tau }(w')}
    \end{aligned}
    \end{equation}
for every probability measure $\bar{\tau }$ on $W^{(0)}$, where $s=|S^{(0)}|$.
    \end{lemm}

    \begin{proof}
By Fubini's theorem,
    \begin{align*}
\int _{\Gamma}&\biggl| \int _{W^{(0)}} \langle a,\pi (M_{\gamma
}w)\rangle \, d\bar{\tau }(w) \biggr|^2 d\gamma
    \\
&= \int_{\Gamma} \int_{W^{(0)}}\int_{W^{(0)}} \langle a,\pi
(M_{\gamma }w)\rangle \overline{\langle a,\pi (M_{\gamma }w')\rangle
}\,d\bar{\tau }(w) \,d\bar{\tau }(w') \,d\gamma
    \\
& = \int_{W^{(0)}}\int_{{W^{(0)}}} \int_{\Gamma} \langle a,\pi
(M_{\gamma } (w-w'))\rangle \,d\gamma \, d\bar{\tau }(w)\, d\bar{\tau
}(w') .
    \end{align*}
 From Lemma \ref{l:character inequality} we know that there exists a
constant $c_a>0$ with
    $$
\int_{\Gamma} \langle a,\pi (M_{\gamma } (w-w'))\rangle \,d\gamma\le
c_a\min(1,\|w-w'\|^{-1/2s})
    $$
for every $w \neq w'$ in $W^{(0)}$, and by integrating we obtain
\eqref{eq:harmonic estimate}.
    \end{proof}


\begin{coro}
\label{corollary about leaf measures}
Let $\bar \tau$ be a probability measure on $W^{ ( 0 ) }$, $x_0\in
X$, and let $\rho = ( \bar \tau
\pi^{-1} ) T_{ - x_0 }$ \textup{(}so $\rho$ is supported on
the central leaf through $x_0$\textup{)}. For every $N \in \N$ we set
\begin{equation*}
\rho_N = \frac 1 N \sum_{ i = 0 }^{N-1} \rho \alpha^i
\end{equation*}
Then for every nontrivial character $a \in \hat X$
\begin{equation*}
\limsup_{ N \to \infty } \absolute { \int \langle a , x \rangle \,d
\rho_N ( x ) }^2 \leq c_a \cdot \int \min ( 1 , \norm { w - w ' }^{ - 1
/ 2 s } ) \,d \bar \tau ( w ) \,d \bar \tau ( w ' ) ,
\end{equation*}
where $c_a$ is as in Lemma \ref{l:character}.
\end{coro}

\begin{proof}
By Cauchy-Schwarz,
\begin{align*}
\absolute { \int \langle a , x \rangle \,d \rho_N ( x ) }^2 & =
\absolute { \frac 1 N \sum_{ i = 0 }^{ N -1} \int \langle a , x \rangle \,d
\rho \alpha^i }^2 \leq \frac{ 1 }{ N } \sum_{ i = 0 }^{N-1} \absolute
{ \int \langle a , x
\rangle \,d \rho a^i ( x ) }^2 \\
& = \frac{ 1 }{ N } \sum_{ i = 0 }^{N-1} \absolute { \int_{ W^{ ( 0 ) } }
\langle a , \pi ( w ) \rangle \,d \bar \tau \bar \beta^i ( w ) }^2
\end{align*}
The map
\begin{equation*}
\gamma \mapsto \absolute { \int_{ W^{ ( 0 ) } } \langle a , \pi ( w )
\rangle \,d ( \bar \tau  M_\gamma ) ( w ) }
\end{equation*}
from $\Gamma$ to $\R^{ + }$ is continuous and bounded, so by the unique
ergodicity of the action of $\bar \beta$ on $\Gamma$
\begin{equation*}
\frac{ 1 }{ N } \sum_{ i = 0 }^{N-1} \absolute { \int_{ W^{ ( 0 ) } }
\langle a , \pi ( w ) \rangle \,d \bar \tau \bar \beta^i ( w ) }^2 \to
\int_\Gamma \absolute { \int_{ W^{ ( 0 ) } } \langle a , \pi ( w )
\rangle \,d \bar \tau M_\gamma }^2 \,d \gamma.
\end{equation*}
We can now apply Lemma \ref{l:character} to conclude the proof of this
corollary.
\end{proof}

    \begin{proof}[Proof of Theorem \ref{t:finite}]
Consider the $\alpha $-invariant Borel set
    $$
B=\{x: \rho _x (W^{(0)}) =\infty\}.
    $$
We will show that
    \begin{equation}
    \label{eq:mu'}
\smash{\mu '=\frac{1}{\mu (B)}\mu\hspace{-1.5mm}\bigm|_B =\lambda _X}
    \end{equation}
whenever $\mu (B)>0$.

Assume therefore that $\mu (B)>0$, and let $X'\subset X$ be the
$\alpha $-invariant Borel set of full measure described in
Proposition \ref{p:rho}. From \eqref{eq:mutranslates} it follows
that, if $x\in B\cap X'$, then any other point in $X'\cap
(x-\pi(W^{(0)}))$ also lies in $B\cap X'$. Hence
    $$
\mu _x^{\mathscr{A}^{(k)}}(B) \in\{0,1\}
    $$
for $\mu \textsl{-a.e.}\;x\in X$ and every $k\ge1$. We can thus
choose, for every $k\ge1$, a set $B^{(k)}\in \mathscr{A}^{(k)}$ with
    $$
\mu (B^{(k)}\bigtriangleup B)=0.
    $$

We fix temporarily a large number $r>0$ and a small $\varepsilon >0$.
According to \eqref{eq:bigballs} there exist an increasing sequence
$(n_k,\,k\ge1)$ of natural numbers and a Borel set $D\subset B$ with
$\mu '(D)> 1-\varepsilon $ so that, for any $x\in D$ and $k\ge1$,
    \begin{gather*}
[x]_{\mathscr{A}^{(n_k)}} +\pi(B_{W^{(0)}}(0,r)) \subset
[x]_{\mathscr{A}^{(n_{k+1})}},
    \\
0<\mu _x^{\mathscr{A}^{(n_{k+1})}}([x]_{\mathscr{A}^{(n_k)}}
+\pi(B_{W^{(0)}}(0,r)) )<\varepsilon
    \end{gather*}
(in the second of these conditions we use the fact that $\rho
_x(W^{(0)})=\infty $ for every $x\in B$). For every $K\ge1$ and $x\in
X$ we set
    $$
\tau_x^K = \frac{1}{K}\sum_{k=1}^K \mu _x^{\mathscr{A}^{(n_k)}}.
    $$
Since $B^{(n_k)}\in\mathscr{A}^{(n_k)}$ is equal to $B$
$(\textup{mod}\;\mu )$ and hence also $(\textup{mod}\;\mu ')$, and
since $\mu '$ is $\alpha $-invariant, we have that
    $$
\mu ' = \int (\mu _x^{\mathscr{A}^{(n_k)}}\alpha ^n)\,d\mu' (x)
    $$
for every $k\ge1$ and $n\in\mathbb{Z}$, and hence that
    \begin{equation}
    \label{eq:integrals}
\mu' =\int (\tau _x^K\alpha ^n) d\mu' (x)
    \end{equation}
for every $K\ge1$ and $n\in\mathbb{Z}$. We define $\bar{\tau }_x^K\in
M_{\infty} (W^{(0)})$ by
    $$
\bar\tau _x^K(C)=\tau _x^K (T_{x}\circ \pi (C))
    $$
for every Borel set $C\subset W^{(0)}$, where $\pi $ and $T_x$ are
taken from \eqref{eq:pi} and \eqref{eq:T}. For any $x\in D$,
    \begin{align*}
M(x,r)&=(\bar\tau_x^K\times\bar\tau_x^K)\bigl( \bigl\{ (w_1,w_2)\in
(W^{(0)})^2: \|w_1-w_2\|<r \bigr\} \bigr)
    \\
&= \smash[b]{\frac{1}{K^2} \cdot \sum_{k=1}^K\sum_{k'=1}^K (\mu
_x^{\mathscr{A}^{(n_k)}}\times\mu _x^{\mathscr{A}^{(n_{k'})}})\bigl(
\bigl\{ (x+\pi(w_1),x+\pi(w_2)):}
    \\
&\qquad\qquad \qquad\qquad \qquad (w_1,w_1)\in (W^{(0)})^2\enspace
\textup{and}\enspace \|w_1-w_2\|<r \bigr\} \bigr)
    \\
&\le \frac{1}{K^2}\cdot \sum_{k=1}^K \;(\mu _x^{\mathscr{A}^{(n_k)}}
\times \mu _x^{\mathscr{A}^{(n_k)}})  (X\times X)
    \\
&\qquad +\smash[b]{\frac{2}{K^2} \cdot \sum_{1\le k<k'\le K} (\mu
_x^{\mathscr{A}^{(n_k)}}\times\mu _x^{\mathscr{A}^{(n_{k'})}}) \bigl(
\bigl\{ (y,y+\pi(w)):}
    \\
&\qquad\qquad \qquad\qquad \qquad\qquad
y\in[x]_{\mathscr{A}^{(n_k)}}, \|w\|<r \bigr\} \bigr)
    \\
&\leq \frac{1}{K} +\frac{2}{K^2}\cdot \sum_{k<k'} \mu
_x^{\mathscr{A}^{(n_{k'})}} \bigl([x]_{\mathscr{A}^{(n_k)}}
+\pi(B_{W^{(0)}}(0,r)) \bigr)<\frac{1}{K}+\varepsilon .
    \end{align*}
Using this, we see that for any $x \in D$,
\begin{equation}
\label{eq:energy estimate}
\begin {aligned}
\int_{ W^{ ( 0 ) } }\int_{ W^{ ( 0 ) } } \min ( 1 &, \norm { w - w '
}^{ - 1 / 2 s } ) \,d
\bar \tau^K_x ( w ) \,d \bar \tau^K_x ( w ' )
    \\
& \leq
M ( x , r ) + r^{ - 1 / 2 s } ( 1 - M ( x , r ) ) \\
& \leq K^{-1} + \varepsilon + r^{ - 1 / 2 s }.
\end {aligned}
\end{equation}

Equation \eqref{eq:integrals} shows that, for an arbitrary positive
integer $K$,
\begin{equation}
\label{eq:intermediate}
\begin {aligned}
\absolute { \int_X \langle a , x \rangle \,d \mu ' }^2 & = \lim_{ N \to
\infty }
\absolute { \int_X \int_X \langle a , y \rangle \,d \left [ \frac{ 1 }{ N
} \sum_{ n = 0 }^{N-1}\bar \tau_x^K \alpha^n \right ] \! ( y ) \, d \mu '
( x ) }^2 \\
& \leq
\int_X \limsup_{ N \to \infty } \absolute { \int_X \langle a , y
\rangle \,d \left [ \frac{ 1 }{ N } \sum_{ n = 0 }^{N-1}\bar \tau_x^K
\alpha^n \right ] \! ( y ) }^2 d \mu ' ( x ) ,
\end {aligned}
\end{equation}
where the first limit is actually over a constant sequence. By
Corollary \ref{corollary about leaf measures},
\eqref{eq:intermediate} and
\eqref{eq:energy estimate},
    \begin{align*}
\biggl| \int_X \langle a,x\rangle \, d\mu'(x) \biggr|&\le
\int_{X\smallsetminus D} \limsup_{ N \to \infty } \absolute { \int_X
\langle a , y
\rangle \,d \left [ \frac{ 1 }{ N } \sum_{ n = 0 }^{N-1} \bar \tau_x^K
\alpha^n \right ] \! ( y ) }^2 d \mu ' ( x ) \\
&\quad +
\int_D\limsup_{ N \to \infty } \absolute { \int_X \langle a , y
\rangle \,d \left [ \frac{ 1 }{ N } \sum_{ n = 0 }^{N-1} \bar \tau_x^K
\alpha^n \right ] \! ( y ) }^2 d \mu ' ( x )
    \\
&\le \mu '(X\smallsetminus D) + c_a\cdot \bigl(K^{-1} +\varepsilon  +
r^{-1/2s}\bigr)
    \\
&\leq \varepsilon  + c_a\cdot \bigl(K^{-1} +\varepsilon  + r^{-1/2s}\bigr).
    \end{align*}
Since $\varepsilon , K,r$ were arbitrary we see that $\int_X \langle
a,x\rangle \, d\mu'(x) =0$ for every $a\in \hat{X}$, and that $\mu'$
is therefore equal to $\lambda _X$.
    \end{proof}

\section{Virtually hyperbolic measures and central equivalence}
    \label{s:virtually hyperbolic}

In this section, we deduce Theorem \ref{t:measures} from Theorem
\ref{t:finite}. For any locally compact metric space $Y$, we let
$M_f(Y) \subset M_\infty (Y)$ denote the finite Borel measures on $Y$.

    \begin{lemm}
    \label{l:center of mass}
There is a Borel map $\mathsf{c}_m : M_f ( \R^d ) \longrightarrow
\R^d$ which commutes with the action of the isometry group of $\R^d$
and is invariant under scalar multiplication: that is, if $F : \R^d
\longrightarrow \R^d$ is an isometry of $\R^d$ and $t > 0$, then
    \begin{equation}
    \label{eq:isometries}
\mathsf{c}_m(\rho ) = F\circ\mathsf{c}_m(t\rho F).
    \end{equation}
    \end{lemm}

    \begin{rema}
    \label{r:center}
For measures $\rho \in M_f(\mathbb{R}^d)$ which have finite first
moments, the vector of moments
    $$
\mathsf{c}_m(\rho )=\biggl(\int x_1\,d\rho (x_1,\dots ,x_d),\dots
,\int x_d\,d\rho (x_1,\dots ,x_d)\biggr)\in\mathbb{R}^d
    $$
would satisfy all these requirements. Unfortunately, there are
measures for which this naive definition of center of mass does not
make sense; the lemma should be interpreted as an alternative,
generalized notion of a center of mass which works for any measure in
$M_f(\mathbb{R}^d)$.
    \end{rema}

    \begin{proof}
For a given $r,\varepsilon >0$, and for every $\rho \in M_f(\mathbb{R}^d)$, let
    \begin{gather*}
S_{r,\varepsilon }(\rho )=\{x\in \mathbb{R}^d:\rho (B(x,r))\ge\varepsilon \},
    \\
a_{r,\varepsilon }(\rho )=\int_{S_{r,\varepsilon }(\rho )}x \,d\rho
(x),\qquad m_{r,\varepsilon }(\rho )=\rho (S_{r,\varepsilon }(\rho )).
    \end{gather*}
Note that the maps $\rho \mapsto m_{r,\varepsilon }(\rho )$ and $\rho
\mapsto \frac{a_{r,\varepsilon }(\rho )}{m_{r,\varepsilon }(\rho )}$,
the latter when defined, are invariant under the action of isometry
group of $\mathbb{R}^d$ on $M_f(\mathbb{R}^d)$. In order to get a map
$\mathsf{c}_m$ which is defined everywhere we arbitrarily fix $r>0$,
set
    $$
n_r(\rho )=\min\,\{n:m_{r,1/n}(\rho )>0\}, \enspace  a_r(\rho
)=a_{r,1/n_r}(\rho ),\enspace  m_r(\rho )=m_{r,1/n_r}(\rho ),
    $$
and put $\mathsf{c}_m(\rho )=\frac{a_r(\rho )}{m_r(\rho )}$.
    \end{proof}

    \begin{proof}
[Proof of Theorem \ref{t:measures}] We first show that for every
$\alpha$-invariant measure $\mu $ which is singular with respect to
the Haar measure $\lambda _X$ there is a virtually hyperbolic measure
$\mu '$ which is centrally equivalent to it. Indeed, consider the map
$\tau : X \longrightarrow X$ defined by
    $$
\tau ( x ) = \pi \circ \mathsf{c}_m ( \rho _x ) + x.
    $$
Let $X '$ be the subset of full measure of $X$ in Proposition
\ref{p:rho}. Then for any $x \in X '$ we have that
    $$
\tau \circ \alpha ( x ) = \pi \circ \mathsf{c}_m ( \mu_{ \alpha x } )
+ \alpha x = \alpha ( \pi \circ \mathsf{c}_m ( \rho _x ) ) + \alpha x
= \alpha \circ \tau ( x )
    $$
where the second equality follows from \eqref{eq:rhoalphabeta} and
\eqref{eq:isometries}. Similarly, by \eqref{eq:mutranslates},
    \begin{equation}
    \label{eq:tau translates}
\begin {aligned} \tau ( x ) & = \pi \circ \mathsf{c}_m ( \rho _x ) +
x=\pi \circ \mathsf{c}_m ( \rho _{y-\pi (w)} ) + x
    \\
& =\pi \circ \mathsf{c}_m (e^{K_\mu (y,w)}\cdot \rho _{y-\pi (w)}) + x
    \\
&=\pi \circ \mathsf{c}_m (\rho _y\bar T_{w}) + x=\pi \circ
\mathsf{c}_m (\rho _y ) + x + \pi ( w ) = \tau (y) \end {aligned}
    \end{equation}
for any $x,y \in X '$ with $y - x = \pi (w) \in X^{(0)}$. By setting
$\mu ' = \mu \tau^{-1}$ we get a new $\alpha $-invariant probability
measure on $X$ which is clearly centrally equivalent to $\mu$. The
$\alpha $-invariant set $Y = \tau (X') \subset X$ is analytic and has
full $\mu '$-measure, and \eqref{eq:tau translates} implies that $Y$
intersects each central leaf in at most one point. By choosing an
$\alpha $-invariant Borel subset $Z\subset Y$ with $\mu '(Z)=1$ we
see that $\mu '$ is virtually hyperbolic.

We now specialize to the case where $\mu $ is ergodic. The map $\tau
^*\colon X'\longrightarrow W^{(0)}$ defined by $\pi \circ \tau
^*(x)=\mathsf{c}_m(\rho _x)$ satisfies that $\tau ^*(\alpha
x)=\bar{\beta }\tau ^*(x)\subset \Gamma \tau ^*(x)$ for every $x\in
X'$ (where $\Gamma $ is defined in \eqref{eq:Gamma}), and the
ergodicity of $\mu $ and the compactness of $\Gamma $ together imply
that there exist an element $w^*\in W^{(0)}$ and a Borel subset
$X''\subset X$ with $\mu (X'')=1$ and $\tau ^*(x)\in\Gamma
w^*=\{M_\gamma w^*:\gamma \in\Gamma \}$ for every $x\in X''$. We
define the probability measure $\tilde{\lambda }_x$ on $X$ as in the
statement of Theorem \ref{t:measures} (3) and obtain that $\mu $ must
be an ergodic component of $\mu ' * \tilde{\lambda }_{\pi (w^*)}$.

For every $\gamma \in\Gamma $ we write $M_\gamma w^*$ as $M_\gamma
w^*=(\gamma _vw^*_v,\,v\in S)$ with $w^*_v=0$ for every $v\in
S\smallsetminus S^{(0)}$ (cf. \eqref{eq:Sc0}). similarly we set $\tau
^*(x)=(\tau ^*(x)_v)$ for every $x\in X''$. Then there exists, for
every $v\in S^{(0)}$ with $w^*_v\ne 0$, a well-defined map $f_v\colon
X'\longrightarrow \mathbb{C}$ with
    $$
\tau ^*(x)_v=f_v(x)w^*_v
    $$
for every $x\in X''$, where we are identifying $K_v$ with
$\mathbb{C}$ (cf. \eqref{eq:iotav}). Since $f_v$ is obviously an
eigenfunction of $\alpha $ for every $v\in S^{(0)}$ with $w^*_v\ne0$,
we have arrived at the following alternative: either the map
$x\mapsto \pi \circ \mathsf{c}_m(x) = \tau ( x )-x$ is zero almost
everywhere, which implies that $\mu = \mu '$, hence virtually
hyperbolic, or $\mu $ is not weakly mixing; indeed, this argument
shows that the point spectrum of $\mu $ (more precisely: the point
spectrum of the action of $\alpha $ on $L^2 (X,\mathscr{S}, \mu )$)
contains some eigenvalue of $\alpha $ of absolute value 1.
    \end{proof}

\section{Central leaves and closed invariant subsets}\label{s:topological}

This section is devoted  to proving the following topological analogue
to Theorem \ref{t:finite}.

\begin{theo}
\label{theorem: topological}
Let $\alpha$ be a nonexpansive, ergodic and totally irreducible
automorphism of a compact connected abelian group $X$. Then any closed
$\alpha$-invariant subset $Y \subsetneq X$ intersects every central
leaf in a compact subset of the leaf.
\end{theo}

The key to this theorem is the following lemma in which we call a
subset $A\subset W^{(0)}$ \emph{$R$-separated} if $\|x-y\|\ge R$ for
any two distinct elements $x,y\in A$.

\begin{lemm}
\label{Lemma: topological}
Let $\alpha$ and $X$ be as in Theorem \ref{theorem: topological}. Then
for any $\varepsilon > 0$ there exist positive integers $R , K$ so
that for any $R$-separated subset
$A \subset W^{ ( 0 ) }$ with at least $K$ elements, the set
\begin{equation*}
\tilde A = \bigcup_{ n = 1 }^\infty \alpha^{ - n } ( \pi ( A ) + x_0 )
\end{equation*}
is $\varepsilon$-dense in $X$.
\end{lemm}

\begin{proof}
Let $\left\{ f_1 , \dots , f_k \right\}$ be a partition of unity of $X$
(i.e. a set of nonnegative continuous functions so that $\sum_{i=1}^k f_i
\equiv 1$) so that the support of each $f_i$ has diameter at most
$\varepsilon$. Clearly, to show that $\tilde A$ is $\varepsilon$-dense it is
sufficient to find some probability measure $\rho$ supported on $\tilde
A$ so that
\begin{equation*}
\int_X f_id \rho > 0
\end{equation*}
for every $i=1,\dots ,k$.
Since the linear span of $\hat X$ is dense in $C ( X )$, there exists
a finite subset $\Xi \subset \hat X$ containing the identity element
$0\in \hat{X}$ so that for each
$i$ we can find an approximation
\begin{equation*}
\tilde f_i ( x ) = \sum_{ a \in \Xi } u_{ i , a } \langle a , x \rangle
\end{equation*}
to $f_i$ in the linear span of $\Xi$ so that
\begin{equation*}
\| f_i - \tilde f_i \|_\infty < \| f_i \|_1 / 100.
\end{equation*}
Let $\Xi ' = \Xi \setminus \{ 0 \}$.

We denote by $c_a$ the constant in Lemma \ref{l:character} and Corollary
\ref{corollary about leaf measures} and define $R$, $K$ by
\begin{equation*}
R^{ 2 s } = K = 100 \max_{ i } \left ( \frac{
\sum_{ a \in \Xi } \absolute { u_{ i , a } c_a } }{ \norm { f_i }_1 }
\right ).
\end{equation*}

Now suppose that $A \subset W^{ ( 0 ) }$ is an $R$-separated set of
cardinality $\geq
K$ and $x_0 \in X$ is arbitrary. We define
\begin{gather*}
\bar \tau = \frac{ 1 }{ \absolute { A } } \sum_{ w \in A } \delta_w ,\enspace
\rho = ( \bar \tau \pi^{-1} ) T_{ - x_0 }, \enspace
\rho_N = \frac{ 1 }{ N } \sum_{ i = 0 }^{N-1} \rho \alpha^i,
\end{gather*}
where $\delta _w$ is the point-mass at $w$. For every $N$, $\rho_N$
is supported on $\tilde A$, and if $N$ is
large enough, then
\begin{align*}
\absolute { \int \langle a , x \rangle \,d \rho_N ( x ) }^2 & \leq 2 c_a
\iint \min ( 1 , \norm { w - w ' }^{ - 1 / 2 s } ) \,d \bar \tau ( w ) \,d
\bar \tau ( w ' ) \\
& \leq 2 c_a ( K^{-1} + R^{ - 1 / 2 s } ) = 4 c_a K^{-1}
\end{align*}
for every $a \in \Xi$.

For $N$ sufficiently large we obtain that
\begin{align*}
\absolute { \int \tilde f_id \rho_N - \int \tilde f_idx } & \leq
\sum_{ a \in \Xi ' } \absolute { u_{ i , a } \int \langle a , x \rangle
d \rho_N } \\
& \leq 4 K^{-1} \sum_{ a \in \Xi ' } \absolute { u_{ i , a } c_a }
\leq \norm { f_i }_1 / 25.
\end{align*}
But then
\begin{align*}
\absolute { \int f_id \rho_N } \geq \absolute { \int \tilde f_id
\rho_N } - \norm { f_i }_1 / 100
\geq \absolute { \int \tilde f_idx } - \norm { f_i }_1 / 20 \geq
\norm { f_i }_1 / 2 > 0
\end{align*}
for $i=1,\dots ,k$, and we are done.
\end{proof}

\begin{proof} [Proof of Theorem \ref{theorem: topological}]
Suppose that $Y\subsetneq X$ is $\alpha$-invariant and closed, and
that the intersection
of $Y$ with some central leaf $X^{ ( 0 ) } + x_0$ is not compact. Fix a
$w_0 \in \pi^{-1} ( x_0 )$ and take $C = [ \pi^{-1} ( Y ) - w_0 ] \cap
W^{ ( 0 ) }$.

By our assumptions, $C$ is a closed unbounded subset of $W^{ ( 0 )
}$. Let $\varepsilon > 0$ be arbitrary, and let $K , R$ be as in Lemma
\ref{Lemma: topological}. Take $A$ to be a finite $R$-separated
subset of $C$ of
cardinality $\geq K$. Then the set $\tilde A \subset X$ defined in
that lemma is a subset of
$Y$ and is $\varepsilon$-dense. So $Y$ is $\varepsilon$-dense, and since
$\varepsilon$ was arbitrary, $Y = X$.
\end{proof}


    \end{document}